\documentclass[11pt]{article}
\usepackage{float}
\usepackage{amsfonts,amsmath,amssymb,graphicx,color}
\usepackage{amsthm}
\usepackage{algorithm}
\usepackage{algpseudocode}
\usepackage{psfrag}
\usepackage{pdflscape}
\usepackage{multirow} 
\usepackage{enumitem}
\usepackage{framed}
\usepackage{sidecap}
\usepackage{tikz}
\usepackage{subcaption}
\usetikzlibrary{arrows}
\usetikzlibrary{patterns,arrows,decorations.pathreplacing}
\graphicspath{ {./figures/} }
    \usepackage{soul}
\usepackage[normalem]{ulem}
\usepackage{url}
\usepackage[numbers, sort&compress ]{natbib} 

\numberwithin{figure}{section}
\numberwithin{table}{section}

\newcommand{\R}{\mathbb{R}}

\newtheorem{thm}{Theorem}[section]
\newtheorem{Def}[thm]{Definition}

\newtheorem{lem}[thm]{Lemma}
\newtheorem{cor}[thm]{Corollary}
\newtheorem{assum}[thm]{Assumptions}
\newtheorem{cond}[thm]{Condition}
\numberwithin{equation}{section}

\newcommand{\comment}[1]{}

\newcommand{\be}{\begin{equation}}
\newcommand{\ee}{\end{equation}}

\newcommand{\bea}{\begin{eqnarray}}
\newcommand{\eea}{\end{eqnarray}}
\newcommand{\beqa}{\begin{eqnarray}}
\newcommand{\eeqa}{\end{eqnarray}}
\newcommand{\beann}{\begin{eqnarray*}}
\newcommand{\eeann}{\end{eqnarray*}}

\newcommand{\bmat}{\left[ \begin{array}}
\newcommand{\emat}{\end{array} \right]}

\newcommand{\beq}{\begin{equation}}
\newcommand{\eeq}{\end{equation}}

\newcommand{\bproof}{\begin{description} \item[{\it Proof}.] ~ }
\newcommand{\eproof}{\hspace*{\fill}$\Box$\medskip \end{description}}

\newcounter{algo}[section]

\newcounter{prog}[section]

\title{A Noise Tolerant SQP Algorithm for Inequality Constrained Optimization}
\author{       
        Figen Oztoprak\thanks{Industrial Engineering Department, Gebze Technical University, Kocaeli, Turkey.  This author was supported by the Scientific and Technological Research Council of Turkey (TUBITAK) under grant 124M396. }
       \and     
        Richard Byrd \thanks{Computer Science Department, University of Colorado, Boulder, USA } 
       }
\begin{document}

\maketitle

\begin{abstract}
We propose a sequential quadratic programming (SQP) algorithm for inequality-constrained optimization that is robust to the presence of bounded noise in function and derivative evaluations.  We cover the case where constraint evaluations contain noise as well as the objective.  The proposed algorithm is a line search SQP method with relaxations to deal with noise.  We study the effect of noise on the global convergence behavior of the algorithm.  We implement the algorithm with noise-aware quasi-Newton updates, and numerically observe that the algorithm can achieve accuracy proportional to the noise level and problem-dependent parameters, as suggested by the theory.
\end{abstract}

\section{Introduction}
\label{intro}
\setcounter{equation}{0}

The majority of mathematical optimization methods are designed for problems where the objective function and the constraint functions defining the feasible set can be evaluated exactly.  However, in several applications, it is either not possible to perform these function and derivative calculations with high precision or is quite computationally expensive.  When the component functions of an optimization problem involve complex approximate computations, contain uncertain varying parameters or involve errors due to simulations, the evaluations of those functions are likely to contain errors.  The derivatives of these functions, however computed, will also return approximate values.  
The effect of these noisy evaluations on optimization algorithms, either existing or new, is increasingly a subject of interest.  However, there are a limited number of studies in the literature that consider the case where constraint evaluations of the optimization problem involve noise \cite{poljak1978nonlinear, hintermuller:2002, conDFO, berghen2005condor, schittkowski2011robust, hansknecht2024, sun2024trust, curtis2025interior}.  

An approach to noise-tolerant optimization is to preserve the core design of standard optimization algorithms with appropriate modifications to noise-sensitive steps \cite{lou2024design}.  This approach has been followed for developing noise-robust methods for unconstrained optimization \cite{berahas2019derivative,shi2022noise,sun2023trust} and for optimization with noisy equality constraints \cite{conDFO, sun2024trust}, under the assumption of bounded noise in function and derivative evaluations.  Recently, Curtis, Dezfulian and Waechter proposed a method based on \cite{curtis2010interior}, which covers bound constraints in addition to noisy equalities \cite{curtis2025interior}.  Their method applies to the subproblems of a nonlinear interior point algorithm for noisy inequality constrained optimization.  In another recent work, Hansknecht, Kirches and Manns proposed a method based on \cite{byrd2003algorithm}, which applies to the nonsmooth exact penalty reformulation of the nonlinear optimization problem with noisy inequality constraints \cite{hansknecht2024}.

In this paper, we propose a line search sequential quadratic programming (SQP) method for nonlinear optimization with noisy inequality constraints.  To make the algorithm robust to the presence of bounded noise, we relax the line search, and modify the quadratic programming (QP) subproblem by allowing relaxations to its feasible set.  We call this algorithm \emph{noise-tolerant SQP with relaxations}(rSQP).  We do not specify the Hessian choice in our algorithm statement. However, in our implementation of rSQP we use quasi-Newton approximations, and provide robustness to noise by a condition introduced in \cite{shi2022noise}.  

We should note that the QP subproblem relaxation we introduced was inspired by the line search penalty-SQP algorithm in \cite{gabo}, although we later noticed its relevance to the SQP algorithm proposed in \cite{burke1989robust}.  The relaxation subproblem \eqref{eq:lpfeas} is also closely related to the vertical step computation of the Byrd-Omojokun trust region SQP algorithm, which is extended and analyzed for noisy equality constraints by Sun and Nocedal \cite{sun2024trust}.  We should also note that another SQP algorithm based on \cite{burke1989robust} for noisy inequality constrained optimization is recently proposed in \cite{facchinei2023}; however, it is designed to work with stochastic noise (which is not necessarily bounded) and therefore follows a different line of analysis.

The existence of inequality constraints adds complications to the analysis of a noise-tolerant SQP.  Unlike the equality-constrained case, regularity cannot be stated in terms of the minimum singular value of the Jacobian.  In this paper, we introduced statements for regularity (see Condition~\eqref{cond:qpreg}), and the failure of regularity (see Definition~\ref{def:almostsingular}) that work with noisy inequalities.  Our analysis in this paper is more extensive compared to our earlier work in \cite{conDFO}, which covered the analysis of an SQP method for noisy equality-constrained optimization under the assumption of regularity only.  

The next section (Section 2) introduces the proposed algorithm after giving the basic definitions we use and the basic assumptions we make on the given optimization problem.  In Section 3, we present its convergence analysis.  Finally, Section 4 covers the details of the numerical implementation of the algorithm and the results of some preliminary tests.

\section{A Noise Tolerant SQP Algorithm}

We consider the following formulation of the general nonlinear optimization problem
\begin{equation}
\label{eq:nlp}
\min_{x\in\R^n} \ \ f(x) \ \ :  \ \ c(x)\leq 0.
\end{equation}
Here, $c:\R^n\rightarrow\R^m$ is a multi-valued function given with $c(x)=\begin{bmatrix}c^1(x)\\ \vdots\\ c^m(x)\end{bmatrix}$.  We pose standard assumptions on \eqref{eq:nlp} as stated in Assumptions~\ref{assum:nlp}, except the one regarding the noisy evaluations.  In the sequel, the noisy evaluations of the functions $f(x)$, $c(x)$, $\nabla f(x)$ and $\nabla c(x)$ at a given point $x$ are denoted with $\tilde f(x)$, $\tilde c(x)$, $\tilde \nabla f(x)$ and $\tilde \nabla c(x)$, respectively. 

\begin{assum}[Assumptions on Problem~\eqref{eq:nlp}]\label{assum:nlp}
The following hold on the domain of functions $f$ and $c$.
\begin{enumerate}
\item Functions $f:\R^n\rightarrow\R$ and $c^i:\R^n\rightarrow\R, i=1,2,\cdots,m$ are continuous and twice continuously differentiable.
\item $f$ is bounded from below.
\item Only noisy function and derivative evaluations of $f$ and $c$ are available such that
\begin{equation}\label{eq:noisebounds1}
|f(x)-\tilde f(x)|\leq \epsilon_f, \qquad \|\nabla\tilde f(x)-\nabla f(x)\|\leq \epsilon_g,  
\end{equation}
and
\begin{equation}\label{eq:noisebounds2}
\|\tilde c(x)-c(x)\|_\infty\leq \epsilon_c, \qquad \|\nabla \tilde c(x)^T-\nabla c(x)^T\|_\infty\leq \epsilon_J. 
\end{equation}
\item The objective function $f$ and the constraint functions $c^i$ have bounded second order derivatives; in particular,
\begin{equation}\label{eq:lipschitzcon}
\|\nabla^2 f(x)\|\leq L_f, \|\nabla^2 c^i(x)\|\leq L_c^i, \ \ \mbox{and} \quad \max\{L_c^1, L_c^2, \cdots, L_c^m\}\leq L_c.
\end{equation}
\item The objective gradient is bounded in norm so that its niosy evaluations are also bounded; i.e.
\[
\|\nabla \tilde f(x)\|\leq G
\]
for some $G<\infty$.
\end{enumerate}
\end{assum}
 
\subsection{Definitions}
We introduce in this subsection the definitions that we need to describe our algorithm.  Other definitions will be given in Section~\ref{sec:analyse} as needed for the analysis of the algorithm.  In the subsequent sections, the sequence of the iterates of the algorithm is denoted with $\{x_k\}$, $k=0,1,2,\cdots$.  Let $g_k=\nabla f(x_k)$, and $J_k=\nabla c(x_k)$ denote evaluations at $x_k$.  Also let $f_k$ and $c_k$ denote $f(x_k)$ and $c(x_k)$, respectively.  Throughout the paper, $\textbf{1}$ denotes the $m$-dimensional vector of ones, $\|.\|_\infty$ denotes the $\ell_\infty$-norm, and $\|.\|$ denotes the $\ell_2$-norm.  The infinity norm of a matrix is the maximum of the $\ell_1$-norm of its rows.

The basic description of our algorithm is not based on a particular norm.  However, to simplify our analysis, we specify the norm to be the $\ell_\infty$-norm in stating the constraint violation, in defining the merit function, and in forming the relaxation subproblem.  We note that the analysis given here could be revised to work with another norm. 

\paragraph{Constraint violation.}
The (non-noisy) measure of infeasibility / constraint violation is given with 
\begin{equation}
v(x) = \|\max\{c(x),0\}\|_\infty.
\end{equation}  
The noisy counterpart is defined as 
\[
\tilde v(x) = \|\max\{\tilde c(x),0\}\|_\infty.
\]
The linear models of $v(x)$ and $\tilde v(x)$ at a given iterate $x_k$ are 
\[
l_v(d;x_k) = \|\max\{[c_k+J_k^Td],0\}\|_\infty, \quad\mbox{and}\quad \tilde l_v(d;x_k) = \|\max\{[\tilde c_k+\tilde J_k^Td],0\}\|_\infty,
\]
respectively.

\paragraph{Merit function.}
We employ the following exact penalty function as the merit function of the SQP algorithm
\begin{equation}
\tilde \phi_{\pi}(x) = \tilde f(x) + \pi(x) \|\max\{\tilde c(x),0\}\|_\infty.
\end{equation}
Here, $\pi(x)$ is the penalty parameter that depends on $x$.  We let $\pi_k=\pi(x_k)$, and define the linear and quadratic models of the merit function as follows.
\[
 \tilde l_{\pi}(d;x_k) = \tilde f_k+\tilde g_k^Td + \pi_k \tilde l_v(d;x_k),
\]
\[
 \tilde q_{\pi}(d;x_k) = \tilde f_k+\tilde g_k^Td + \frac{1}{2}d^TH_kd + \pi_k \tilde l_v(d;x_k).
\]
Here, $H_k$ is an iteration dependent positive definite matrix specified in the algorithm.  We also note here the non-noisy counterparts of the above models.
\[
 l_{\pi}(d;x_k) = f_k+g_k^Td + \pi_k l_v(d;x_k),
\]
and
\[
 q_{\pi}(d;x_k) = f_k+g_k^Td + \frac{1}{2}d^TH_kd + \pi_k l_v(d;x_k).
\]

\subsection{Algorithm Description }

\begin{algorithm}[h!] 
\caption{Noise Tolerant SQP with Relaxations (rSQP)}
\label{alg:rsqp}

\begin{algorithmic}[1]
\State \textit{[Initialize.]} Choose $\Delta_{\min}<\Delta_0<\Delta_{\max}$, $\pi_0\geq 1$, and $0<\theta_1,\theta_2<1$. Set k=0.

\State \textit{[Compute relaxation.]} Solve
\begin{align}
\label{eq:lpfeas}
(r_v, d_v) = \mbox{argmin}_{(r,d)} \ & \|r\|_\infty \nonumber\\
 s.t. \quad   & \tilde c_k+\tilde J_k^Td\leq r\\
     & r \geq 0, \|d\|\leq \Delta_k.\nonumber
\end{align}

\State \textit{[Compute direction.]} 
For $H_k$ positive definite solve
\begin{align}
\label{eq:qpdirobj}
d_{qp} = \mbox{argmin}_d \ \ & \tilde f_k + \tilde g_k^Td + \frac{1}{2}d^TH_kd\\
\label{eq:qpdircon}
s.t. \quad & \tilde c_k+\tilde J_k^Td \leq r_v.
\end{align}

\State \textit{[Update penalty parameter.]} If $\tilde l_v(0;x_k) - \tilde l_v(d_{qp};x_k) = \tilde v(x_k) - \|r_v\| > 0$, find smallest $\pi\geq 1.1    \pi_k$ satisfying
\begin{equation}\label{piupdate}
\tilde q_{\pi}(0;x_k) - \tilde q_{\pi}(d_{qp};x_k)] \geq \theta_1\pi[\tilde l_v(0;x_k) - \tilde l_v(d_{qp};x_k)].
\end{equation}
Set $\pi_k=\pi$.

\State \textit{[Line search.]} Set $\alpha = 2^{-j}$, where $j$ is the smallest nonnegative integer such that
\begin{equation}\label{eq:lscond}
 \tilde\phi_{\pi}(x_k) - \tilde\phi_{\pi}(x_k+\alpha d_{qp}) \geq \theta_2\alpha[\tilde q_{\pi}(0;x_k) - \tilde q_{\pi}(d_{qp};x_k)] - 2\epsilon_R,
\end{equation}
where $\epsilon_R = \epsilon_f+\pi_k\epsilon_c$.  

\bigskip

\State Set $x_{k+1}=x_k+\alpha_k d_{qp}$.  $k = k+1$.  Choose $\Delta_{\min}\leq\Delta_k\leq\Delta_{\max}$.  Go to Step 2.
\end{algorithmic} 
\end{algorithm}

The proposed noise-tolerant algorithm is an extension of conventional line search SQP.  To tolerate noise, we make two modifications to the base algorithm.  First, the same way as in \cite{conDFO}, we relax the line search in a way that it does not fail due to noise.  Second, we relax the feasible region of the QP subproblem such that it is feasible regardless of the effect of the noise in evalutions.    

A complete description of the algorithm is given in Algorithm~\ref{alg:rsqp}.  The search direction is obtained via solving a convex QP in Step 3.  The consistency of the constraints is guaranteed by introducing a relaxation term $r_v$.  The amount of relaxation is computed in Step 2 via minimizing the linear model of the constraint violation given the noisy evaluations $\tilde c_k$, and $\tilde J_k$.  

Note that any optimal solution to \eqref{eq:lpfeas} is feasible for the QP subproblem \eqref{eq:qpdirobj}-\eqref{eq:qpdircon}, and is bounded in norm by $\Delta_{\max}$.  The construction of \eqref{eq:lpfeas} is indeed the same as the linear programming subproblem of the penalty-SQP algorithm in \cite{gabo}.  However, we use it for directly relaxing the feasible set of the QP subproblem rather than for setting a penalty parameter value to provide that.

\section{Convergence Analysis}\label{sec:analyse}

We analyze the convergence behavior of Algorithm~\ref{alg:rsqp} under Assumptions~\ref{assum:nlp} regarding problem \eqref{eq:nlp}, and the following assumptions on the algorithm.  Basically, we require that the QP subproblems have strictly convex objectives.  We also assume for simplicity that the parameter $\Delta_k$ of subproblem \eqref{eq:lpfeas} is fixed to a positive value for all $k$.   
\begin{assum}[Assumptions on Algorithm~\ref{alg:rsqp}]\label{assum:alg}
The following hold for Algorithm~\ref{alg:rsqp}.
\begin{enumerate}
\item $H_k$ satisfies $\sigma_1 I \succeq H_k \succeq \sigma_n I$, $\sigma_1\geq\sigma_n> 0$, $\forall k$.
\item $\Delta_{\min}=\Delta_{\max}=\Delta$ for some $\Delta>0$.
\item $\pi_k\geq 1$, $\forall k$.
\end{enumerate}
\end{assum}

\subsection{Boundedness of $d_{qp}$}

We first present some preliminary results on the direction $d_{qp}$ computed in Step 3 of Algorithm~\ref{alg:rsqp}.  

\begin{lem}\label{lem:welldqp}
Suppose Assumptions~\ref{assum:nlp} and \ref{assum:alg} hold.  A solution $d_{qp}$ to subproblem \eqref{eq:qpdirobj}-\eqref{eq:qpdircon} exists and is unique at all iterations of Algorithm~\ref{alg:rsqp}.
\end{lem}

The proof is immediate as \eqref{eq:qpdirobj}-\eqref{eq:qpdircon} is always feasible by construction, and is a convex program with a strictly convex objective by Assumptions~\ref{assum:alg}.  

\bigskip

Note that subproblem \eqref{eq:lpfeas} is equivalent to 
\[
\min_{\|d\|\leq\Delta} \  \tilde l_v(d;x_k) =  \|\max\{[\tilde c_k+\tilde J_k^Td],0\}\|_\infty.
\]
By \eqref{eq:qpdircon}, the linear model improvement in constraint violation provided by $d_{qp}$ cannot be less than the improvement by any solution $d_v$ to \eqref{eq:lpfeas}.  We state that formally in Lemma~\ref{lem:feasdqp}.  

\begin{lem}\label{lem:feasdqp}
At all iterations of Algorithm~\ref{alg:rsqp}, the computed direction $d_{qp}$ satisfies
\[
\tilde v(x_k) - \tilde l_v(d_{qp};x_k) \geq \tilde v(x_k) - \tilde l_v(d_v;x_k).
\]
\end{lem}

The proof follows directly from the formulations of \eqref{eq:lpfeas} and \eqref{eq:qpdircon}.  The next lemma provides a bound on the norm of the direction $d_{qp}$.

\begin{lem}\label{lem:lendqp}
Suppose Assumptions~\ref{assum:nlp} and \ref{assum:alg} hold.  Then, 
\[
\|d_{qp}\| \leq \frac{4G+\sigma_1\Delta}{\sigma_n} =:r_l.
\]
\end{lem}
\bproof
Let $\bar d$ denote the minimum norm solution to \eqref{eq:lpfeas}; by construction, $\|\bar d\|\leq \Delta$.  Following Lemma~\ref{lem:feasdqp}, we observe that $\|d_{qp}\|\geq\|\bar d\|$ because otherwise $\bar d$ would not be the minimum norm solution.  On the other hand, by the optimality of $d_{qp}$,
\begin{align*}
& \tilde g_k^Td_{qp} + \frac{1}{2}d_{qp}^TH_kd_{qp} \leq \tilde g_k^T\bar d + \frac{1}{2}\bar d^TH_k\bar d \\
\Rightarrow  \ & \frac{1}{2}\sigma_n\|d_{qp}\|^2\leq  \frac{1}{2}d_{qp}^TH_kd_{qp} \leq \tilde g_k^T(\bar d -d_{qp})+ \frac{1}{2}\bar d^TH_k\bar d \leq \tilde g_k^T(\bar d -d_{qp})+ \frac{1}{2}\sigma_1\|\bar d\|^2\\
\Rightarrow  \ & \frac{1}{2}\sigma_n\|d_{qp}\|^2\leq G(\|d_{qp}\|+\|\bar d\|)+ \frac{1}{2}\sigma_1\|\bar d\|^2\\
\Rightarrow  \ & \frac{1}{2}\sigma_n\|d_{qp}\|^2 - 2G\|d_{qp}\|-\frac{1}{2}\sigma_1\Delta^2\leq 0.
\end{align*}
The resulting convex quadratic inequality is satisfied for  
\[
\|d_{qp}\|\leq \frac{2G+\sqrt{4G^2+\sigma_n\sigma_1\Delta^2}}{\sigma_n}, 
\]
which implies 
\[
\|d_{qp}\| \leq \frac{4G+\sigma_1\Delta}{\sigma_n}.
\]
\eproof

\subsection{Feasibility}\label{sec:feas}

\paragraph{Stationarity measure for minimization of constraint violation.}  To cover the case when the algorithm converges to an infeasible stationary point, we need to consider a stationarity measure for the minimization of $v(x)$.  Considering a fixed $\Delta>0$ we define 
\[
\psi_v(x_k) = v(x_k)-\min_{\|d\|\leq \Delta} \ l_v(d;x_k).
\] 
Also define the corresponding noisy measure as follows.
\[
\tilde\psi_v(x_k) = \tilde v(x_k)-\min_{\|d\|\leq \Delta} \ \tilde l_v(d;x_k).
\] 
Observe that Lemma~\ref{lem:feasdqp} links $d_{qp}$ to the noisy measure $\tilde\psi_v(x_k)$.  In particular, this lemma states that the reduction achieved by $d_{qp}$ in the linear model of constraint violation is not less than the value of the noisy stationarity measure; i.e.,  
\[
\tilde\psi_v(x_k) = \tilde v(x_k) - \tilde l_v(d_v;x_k) \leq \tilde v(x_k) - \tilde l_v(d_{qp};x_k).
\]
The next result establishes a bound on the difference of $\tilde\psi_v(x_k)$ and $\psi_v(x_k)$.  We readily know that 
\[
0\leq \psi_v(x_k) \leq v(x_k),
\]
and
\[
0\leq \tilde\psi_v(x_k) \leq \tilde v(x_k).
\]
Consider the extreme case where $\tilde\psi_v(x_k)=0$ and $\psi_v(x_k)=v(x_k)$ (which seems possible for some $\epsilon_J$).  Or, $l_v(x;x_k)=0$ for some $x$ satisfying $\|x-x_k\|\leq \Delta$, whereas $\tilde l_v(x;x_k)\geq \tilde v(x_k)$ for all $x$.  Lemma~\ref{lem:diffpsiv} implies that this can occur only if $\Delta$ and $\epsilon_J$ are large enough.
\begin{lem}\label{lem:diffpsiv}
Suppose Assumptions~\ref{assum:nlp} hold.  Let $E_\epsilon^v = 2\epsilon_c+\Delta\epsilon_J$.  Then, 
\[
|\tilde \psi_v(x_k) -\psi_v(x_k) | \leq E_\epsilon^v.
\] 
\end{lem}
\bproof
Let 
$d_v \in \mbox{argmin}_{\|d\|\leq \Delta} \ l_v(d;x_k)$,
and 
$\tilde d_v \in \mbox{argmin}_{\|d\|\leq \Delta} \ \tilde l_v(d;x_k)$.
Since both $\|\tilde d_v\|\leq\Delta$ and $\|d_v\|\leq\Delta$ hold, we have
\begin{align*}
\tilde l_v(\tilde d_v;x_k) - l_v(d_v;x_k) & = \tilde l_v(\tilde d_v;x_k) -\tilde l_v(d_v;x_k) +\tilde l_v(d_v;x_k) - l_v(d_v;x_k)\\
& \leq \tilde l_v(d_v;x_k) - l_v(d_v;x_k)\\ 
& = \|\max\{\tilde c_k+\tilde J_k^Td_v,0\}\|_\infty-\|\max\{c_k+J_k^Td_v,0\}\|_\infty \\
& \leq \|\max\{\tilde c_k+\tilde J_k^Td_v,0\}-\max\{c_k+J_k^Td_v,0\}\|_\infty \\
& \leq \|(\tilde c_k+\tilde J_k^Td_v)-(c_k+J_k^Td_v)\|_\infty \\
& \leq \|\tilde c_k-c_k\|_\infty+\|\tilde J_k^T-J_k^T\|_\infty\|d_v\|_\infty\\
& \leq \epsilon_c+\epsilon_J\|d_v\| \leq \epsilon_c+\Delta\epsilon_J.
\end{align*}
The same bound holds for $l_v(d_v;x_k)-\tilde l_v(\tilde d_v;x_k)$.  Therefore,
\begin{equation}\label{eq:lvbound}
|\tilde l_v(\tilde d_v;x_k) - l_v(d_v;x_k) | \leq \epsilon_c+\Delta\epsilon_J,
\end{equation}
and
\[
|\tilde \psi_v(x_k) -\psi_v(x_k) | \leq |\tilde v(x_k)-\tilde l_v(\tilde d_v;x_k) - (v(x_k) - l_v(d_v;x_k))| \leq 2\epsilon_c+\Delta\epsilon_J.
\]
\eproof

\paragraph{A weak bound on the steplength.}
In our analysis, we will give two bounds on the steplength computed in Step 5 of Algorithm~\ref{alg:rsqp}.  The one that we will give in this section is the weaker bound, which always holds independent of regularity.  We will provide in Section~\ref{sec:opt} a stronger bound that applies only when regularity conditions hold.

\begin{lem}\label{lem:deltabound}
Let $\pi=\pi_k$, and suppose Assumptions \ref{assum:nlp} hold.  Define
\begin{equation} \label{eq:deltabound}
\delta_k = q_\pi(0;x_k) - q_\pi(d_{qp};x_k) - [\tilde q_\pi(0;x_k) - \tilde q_\pi(d_{qp};x_k)]. 
\end{equation}
Then,
\[
\delta_k \geq -\epsilon_g r_l - 2\pi\epsilon_c -\pi\epsilon_Jr_l.
\]
\end{lem}

\bproof
Observe that 
\begin{align*}
v(x_k)-\tilde v(x_k)& =\|\max\{c_k,0\}\|_\infty-\|\max\{\tilde c_k,0\}\|_\infty\\
& \geq -\|\max\{\tilde c_k,0\}-\max\{c_k,0\}\|_\infty\geq-\|\tilde c_k-c_k\|_\infty.
\end{align*}
Similarly,
\[
\tilde l_v(d_{qp};x_k)-l_v(d_{qp};x_k) \geq -\|(\tilde c_k+\tilde J_k^Td_{dqp})-(c_k + J_k^Td_{qp})\|_\infty.
\]
By \eqref{eq:noisebounds1}-\eqref{eq:noisebounds2}, Lemma~\ref{lem:lendqp}, and using the fact that $\|d_{qp}\|_\infty\leq \|d_{qp}\|$ we obtain
\begin{align*}
\delta_k & = -(g_k - \tilde g_k)^Td_{qp} + \pi [v(x_k)-\tilde v(x_k)]+ \pi[\tilde l_v(d_{qp};x_k)-l_v(d_{qp};x_k)]\\
& \geq -\|g_k - \tilde g_k\|\|d_{qp}\| -\pi\|\tilde c_k - c_k\|_\infty - \pi\|(\tilde c_k - c_k)+(\tilde J_k - J_k)^Td_{qp}\|_\infty\\
& \geq -\epsilon_g r_l - 2\pi\epsilon_c -\pi\epsilon_Jr_l.
\end{align*}
\eproof

\begin{lem}\label{lem:alpha_weak}
Let $\pi=\pi_k$ and $L_\pi = L_f+\pi L_c$.  Define $E^v_\epsilon$ as in Lemma~\ref{lem:diffpsiv} and suppose Assumptions \ref{assum:nlp} and \ref{assum:alg} hold.  At any iteration $k$ of Algorithm~\ref{alg:rsqp}, if
\begin{equation}\label{eq:vcond}
\psi_v(x_k)\geq \frac{2}{\tau}\max\left\lbrace E^v_\epsilon,\frac{\theta_1(1-\theta_2)}{\theta_2}\left( 2\epsilon_c+r_l\epsilon_J + \frac{1}{\pi}\epsilon_gr_l\right) \right\rbrace 
\end{equation}
for some $\tau\in(0,1)$, then the line search step of Algorithm~\ref{alg:rsqp} will return 
\begin{equation}
\alpha_k \geq \frac{(1-\tau)(1-\theta_2)\theta_1\psi_v(x_k)}{(L_c + \frac{1}{\pi}L_f)r_l^2} + \frac{\sigma_n}{2L_\pi}.
\end{equation}  
\end{lem}

\bproof
Using the definition of $\delta_k$ in \eqref{eq:deltabound},
\begin{align*}
\phi_{\pi}(x_k) - \phi_{\pi}(x_k+\alpha d_{qp}) & \geq l_\pi(0;x_k) - l_\pi(\alpha d_{qp};x_k) -\frac{L_\pi}{2}\alpha^2\|d_{qp}\|^2\\
& \geq \alpha[l_\pi(0;x_k) - l_\pi(d_{qp};x_k)]-\frac{L_\pi}{2}\alpha^2\|d_{qp}\|^2\\
& = \alpha[q_\pi(0;x_k) - q_\pi(d_{qp};x_k)]+\frac{\alpha}{2}d_{qp}H_kd_{qp}-\frac{L_\pi}{2}\alpha^2\|d_{qp}\|^2\\
& \geq \alpha[q_\pi(0;x_k) - q_\pi(d_{qp};x_k)]+\alpha\frac{\sigma_n}{2} \|d_{qp}\|^2-\frac{L_\pi}{2}\alpha^2 \|d_{qp}\|^2\\
& \geq \alpha[\tilde q_\pi(0;x_k) - \tilde q_\pi(d_{qp};x_k)+\delta_k]+\alpha\frac{\sigma_n}{2} \|d_{qp}\|^2-\frac{L_\pi}{2}\alpha^2\|d_{qp}\|^2\\
& \geq \alpha\theta_2[\tilde q_\pi(0;x_k) - \tilde q_\pi(d_{qp};x_k)]+\alpha\theta_2\delta_k +  \alpha(1-\theta_2)[\tilde q_\pi(0;x_k) - \tilde q_\pi(d_{qp};x_k)]\\
& \qquad +\alpha\frac{\sigma_n}{2} \|d_{qp}\|^2-\frac{L_\pi}{2}\alpha^2\|d_{qp}\|^2\\
& \geq \alpha\theta_2[\tilde q_\pi(0;x_k) - \tilde q_\pi(d_{qp};x_k)]+\alpha\theta_2\delta_k +  \alpha(1-\theta_2)\theta_1\pi[\tilde l_v(0;x_k) - \tilde l_v(d_{qp};x_k)]\\
& \qquad +\alpha\frac{\sigma_n}{2} \|d_{qp}\|^2-\frac{L_\pi}{2}\alpha^2\|d_{qp}\|^2.\\
\end{align*}
The last line follows from \eqref{piupdate}.  On the left side of the inequality, by choice of $\epsilon_R$ and the assumptions \eqref{eq:noisebounds1} and \eqref{eq:noisebounds2} we obtain
\begin{equation}\label{eq:eqrelax}
\tilde\phi_{\pi}(x_k) - \tilde\phi_{\pi}(x_k+\alpha d_{qp}) \geq \phi_{\pi}(x_k) - \phi_{\pi}(x_k+\alpha d_{qp})-2\epsilon_R.
\end{equation}
Therefore, the line search condition \eqref{eq:lscond} is satisfied for $\alpha$ such that
\[
\alpha\theta_2\delta + \alpha(1-\theta_2)\theta_1\pi[\tilde l_v(0;x_k) - \tilde l_v(d_{qp};x_k)] +\alpha\frac{\sigma_n}{2} \|d_{qp}\|^2
- \frac{L_\pi}{2}\alpha^2\|d_{qp}\|^2 \geq 0,
\]
or $\alpha>0$ with
\begin{equation}\label{eq:amid}
\alpha \leq 2\frac{\theta_2\delta + (1-\theta_2)\theta_1\pi[\tilde l_v(0;x_k) - \tilde l_v(d_{qp};x_k)]}{L_\pi\|d_{qp}\|^2} + \frac{\sigma_n}{L_\pi} = 2\frac{\theta_2\delta + (1-\theta_2)\theta_1\pi\tilde \psi_v(x_k)}{L_\pi\|d_{qp}\|^2}+ \frac{\sigma_n}{L_\pi}. 
\end{equation}
Using Lemmas~\ref{lem:diffpsiv} and \ref{lem:deltabound}, we have 
\begin{align*}
\theta_2\delta + (1-\theta_2)\theta_1\pi\tilde \psi_v(x_k) & \geq (1-\theta_2)\theta_1\pi\psi_v(x_k)- (1-\theta_2)\theta_1\pi E_\epsilon^v + \theta_2\delta \\
& \geq (1-\theta_2)\theta_1\pi\psi_v(x_k)-(1-\theta_2)\theta_1\pi E_\epsilon^v -\pi\theta_2\left(\frac{1}{\pi}\epsilon_g r_l +(2\epsilon_c +\epsilon_Jr_l)\right)\\
& \geq (1-\theta_2)\theta_1\pi\psi_v(x_k)-(1-\theta_2)\theta_1\pi\tau\psi_v(x_k).
\end{align*}
As we also have $\|d_{qp}\|\leq r_l$ by Lemma~\ref{lem:lendqp}, we conclude that \eqref{eq:amid} holds for all $\alpha$ satisfying
\[
\alpha \leq 2\frac{(1-\tau)(1-\theta_2)\theta_1\pi\psi_v(x_k)}{L_{\pi} r_l^2}+ \frac{\sigma_n}{L_\pi}=2\frac{(1-\tau)(1-\theta_2)\theta_1\psi_v(x_k)}{(L_c + \frac{1}{\pi}L_f) r_l^2}+ \frac{\sigma_n}{L_\pi}.
\]
The line search step of Algorithm~\ref{alg:rsqp} halves the steplength at each trial, therefore the resulting $\alpha_k$ cannot be smaller than the half of the above upper bound.
\eproof

\paragraph{Critical region for feasibility.}
We are now ready to state a result on the sequence $\{\psi_v(x_k)\}$.  Let us start by introducing 
\[
E_\epsilon^{\mathcal{F}}= \max\left\lbrace E^v_\epsilon,\frac{\theta_1(1-\theta_2)}{\theta_2}\left( 2\epsilon_c+r_l\epsilon_J + \frac{1}{\pi_0}\epsilon_gr_l\right) \right\rbrace.
\]
Here, $E^v_\epsilon=2\epsilon_c+\Delta\epsilon_J$ as given in Lemma~\ref{lem:diffpsiv}.  Note that $E_\epsilon^{\mathcal{F}}$ is an upper bound to the right hand side of \eqref{eq:vcond} since $\pi_k$ is not allowed to decrease. 
\begin{thm}\label{thm:crfeas}
Consider the iterates $\{x_k\}$ of Algorithm~\ref{alg:rsqp}, and suppose Assumptions \ref{assum:nlp} and \ref{assum:alg} hold.  The algorithm must include infinitely many iterates in the region $\mathcal{F}$ (a critical region for feasibility) defined by
\begin{equation}\label{eq:crit1}
\mathcal{F} = \left\{x: \psi_v(x) \leq  \max\left(\frac{2}{\tau}E_\epsilon^{\mathcal{F}}, 2\sqrt{\gamma\left(\frac{1}{\pi_0}\epsilon_f+\epsilon_c\right)} \right)\right\}
\end{equation}
for $\gamma=\displaystyle\frac{2(L_c + \frac{1}{\pi_0}L_f)r_l^2}{\tau(1-\tau)(1-\theta_2)\theta_2\theta_1^2}$ and $\tau\in(0,1)$.
Each time the iterates leave $\mathcal{F}$, they will return  to $\mathcal{F}$.

\end{thm}
\bproof
Consider the following scaled merit function, and its noisy counterpart 
\begin{equation}\label{eq:scaledmerit}
\phi_{\pi}^s(x_k) = \frac{1}{\pi} (f_k+F) + \|\max\{c_k,0\}\|_\infty, 
\qquad \tilde\phi_{\pi}^s(x_k) = \frac{1}{\pi} (\tilde f_k+F) + \|\max\{\tilde c_k,0\}\|_\infty.
\end{equation}
Here, $F$ is a sufficiently large value so that  $f_k+F\geq 0$ and $\tilde f_k+F\geq 0$ for all $k$ (i.e. $\phi_{\pi}^s(x_k)\geq 0$ and $\tilde \phi_{\pi}^s(x_k)\geq 0$).  Observe that 
\[
\phi_{\pi}^s(x_k)-\phi_{\pi}^s(x_k+\alpha d) = \frac{1}{\pi_k} [\phi_{\pi}(x_k)-\phi_{\pi}(x_k+\alpha d)]
\]
so that the following hold for the steplengths $\alpha_k$ obtained by the algorithm:
\begin{equation}\label{eq:intres2}
\tilde\phi^s_{\pi}(x_k) - \tilde\phi^s_{\pi}(x_k+\alpha_k d_{qp}) \geq \theta_2\alpha_k[\tilde q^s_{\pi}(0;x_k) - \tilde q^s_{\pi}(d_{qp};x_k)] - \frac{2}{\pi_k}\epsilon_R.
\end{equation}
Moreover,
\begin{align*}
\tilde\phi^s_{\pi_k}(x_k) - \tilde\phi^s_{\pi_{k+1}}(x_{k+1}) &= \frac{1}{\pi_k} (\tilde f_k+F) + \|\max\{\tilde c_k,0\}\|_\infty - \frac{1}{\pi_{k+1}} (\tilde f_{k+1}+F) - \|\max\{\tilde c_{k+1},0\}\|_\infty\\
&\geq \frac{1}{\pi_k} (\tilde f_k+F) + \|\max\{\tilde c_k,0\}\|_\infty - \frac{1}{\pi_{k}} (\tilde f_{k+1}+F) - \|\max\{\tilde c_{k+1},0\}\|_\infty\\
& = \tilde\phi^s_{\pi_k}(x_k) - \tilde\phi^s_{\pi_k}(x_{k+1})
\end{align*}
as $\pi_k\leq \pi_{k+1}$.  Therefore, the decrease in $\tilde{\phi}^s_{\pi_k}(x_k)$ is monotonic even when $\pi_k$ keeps growing.

\bigskip

Now, note that the update condition \eqref{piupdate} for the penalty parameter implies
\begin{equation}\label{eq:intres1}
\tilde q^s_\pi(0;x_k) - \tilde q^s_\pi(d_{qp};x_k) \geq \theta_1[\tilde l_v(0;x_k) - \tilde l_v(d_{qp};x_k)] = \theta_1\tilde \psi_v(x_k).
\end{equation}
On the other hand, we have 
\[
\phi^s_{\pi}(x_k) - \phi^s_{\pi}(x_k+\alpha_k d_{qp}) +\frac{2}{\pi_k}\epsilon_R \geq \tilde\phi^s_{\pi}(x_k) - \tilde\phi^s_{\pi}(x_k+\alpha_k d_{qp}),
\]
which together with \eqref{eq:intres1}, \eqref{eq:intres2}, Lemma~\ref{lem:diffpsiv}, and Lemma~\ref{lem:alpha_weak} yields
\begin{align*}
& \phi^s_{\pi}(x_k) - \phi^s_{\pi}(x_k+\alpha_k d_{qp}) + \frac{4}{\pi_k}\epsilon_R \geq \theta_1\tilde \psi_v(x_k) \geq \theta_1\theta_2\alpha_k(\psi_v(x_k)-E_\epsilon^v)\\
\Rightarrow \ & \phi^s_{\pi}(x_k) - \phi^s_{\pi}(x_k+\alpha_k d_{qp}) \geq \theta_1\theta_2\alpha_k(\psi_v(x_k)-E_\epsilon^v) - \frac{4}{\pi_k}\epsilon_R\\
& \geq \frac{(1-\tau)(1-\theta_2)\theta_2\theta_1^2\psi_v(x_k)}{(L_c + \frac{1}{\pi_0}L_f)r_l^2}(\psi_v(x_k)-E_\epsilon^v) +\frac{\theta_1\theta_2\sigma_n}{2L_\pi}(\psi_v(x_k)-E_\epsilon^v) - \frac{4}{\pi_k}\epsilon_R
\end{align*}
By \eqref{eq:crit1}, for $x_k\notin \mathcal{F}$ we have
\[
\frac{\tau}{2}\psi_v(x_k) \geq E_\epsilon^v, 
\]
and
\[
\frac{\tau}{2}\frac{(1-\tau)(1-\theta_2)\theta_2\theta_1^2}{(L_c + \frac{1}{\pi_0}L_f)r_l^2}\psi_v(x_k)^2 \geq 4\left(\frac{1}{\pi_0}\epsilon_f+\epsilon_c\right)\geq \frac{4}{\pi_k}\epsilon_R.
\]
So, we end up with
\begin{align*}
\phi^s_{\pi}(x_k) - \phi^s_{\pi}(x_{k+1}) &\geq (1-\tau)^2\frac{(1-\theta_2)\theta_2\theta_1^2}{(L_c + \frac{1}{\pi_0}L_f)r_l^2}\psi_v(x_k)^2 +(1-\frac{\tau}{2})\frac{\theta_1\theta_2\sigma_n}{2\pi_k(L_c + \frac{1}{\pi_0}L_f)}\psi_v(x_k).
\end{align*}
Since $\phi^s_{\pi}$ is bounded from below, if  there exists an index $k'$ such that $x_{k} \notin \mathcal{F}$ for all $k \geq k'$ then  $\psi_v(x_k)$ is forced to converge to zero.  Since this cannot occur, we must have  
 $x_{k} \in \mathcal{F}$ for some iterate $k >k'$ .
 Therefore the algorithm produces infinitely many iterates in $\mathcal{F}$.
\eproof

Note that the statement of $\mathcal{F}$ contains components that are relevant to the objective function.  These components would disappear if the objective function is non-noisy, or they would fade if the penalty parameter tends to infinity.

\bigskip

\subsection{Optimality}\label{sec:opt}
\paragraph{Stationarity measure for optimality.}
Let $\tilde D$ denote the feasible region of the QP \eqref{eq:qpdirobj}-\eqref{eq:qpdircon}; i.e.
\begin{equation}\label{eq:tildeD}
\tilde D = \{d : \tilde c_k+\tilde J_k^Td \leq r_v , \ r_v \in \underset{r\geq \tilde c_k+\tilde J_k^Tv, \ \|v\|\leq \Delta}{\mbox{argmin}} \|r\|_\infty\}.
\end{equation}
Define
\[
\tilde q_f(d;x_k,H_k) = \tilde f_k + \tilde g_k^Td + \frac{1}{2}d^TH_kd,
\]
and
\[
\tilde q_f(d;x_k,\sigma) = \tilde f_k + \tilde g_k^Td + \frac{\sigma}{2}\|d\|^2.
\]
We consider the following two functions to measure distance to optimality.
\begin{align*}
\tilde \psi_{o}(x_k;H_k) & = \tilde f_k - \min_{d\in \tilde D} \ \tilde q_f(d;x_k,H_k) = \tilde g_k^Td_{qp} + \frac{1}{2}d_{qp}^TH_kd_{qp},\\
\tilde \psi_{o}(x_k) & = \tilde f_k - \min_{d\in \tilde D} \ \tilde q_f(d;x_k,\sigma_n).
\end{align*}
We introduce $\tilde\psi_{o}(x_k)$ since $\tilde\psi_{o}(x_k;H_k)$ depends on $H_k$.  On the other hand, we prefer to use $\tilde\psi_{o}(x_k;H_k)$ as it is readily computed by the algorithm.  

Likewise, we introduce the corresponding non-noisy functions
\[
q_f(d;x_k,H_k) = f_k + g_k^Td + \frac{1}{2}d^TH_kd, \quad q_f(d;x_k,\sigma) = f_k + g_k^Td + \frac{\sigma}{2}\|d\|^2,
\]
and
\[
\psi_{o}(x_k;H_k) = f_k - \min_{d\in D} \ q_f(d;x_k,H_k), \quad  \psi_{o}(x_k) = f_k - \min_{d\in D} \ q_f(d;x_k,\sigma_n),
\]
where
\begin{equation}\label{eq:D}
D = \{d : c_k+J_k^Td \leq r_v , \ r_v \in \underset{r \geq c_k+J_k^Tv, \ \|v\|\leq \Delta} {\mbox{argmin}} \|r\|_\infty\}.
\end{equation}

We first discuss in Lemma~\ref{lem:psioval} that the proposed optimality measures are valid.  Next, in Lemma~\ref{lem:psiorel}, we compare $\psi_{o}(x_k;H_k)$ to $\psi_{o}(x_k)$.  Finally, in Lemma~\ref{lem:psiodiff}, we give an upper bound on the difference in between the noisy and non-noisy measures computed at a given point. 

\begin{lem}\label{lem:psioval}
Suppose $\bar x$ is a feasible solution to \eqref{eq:nlp} where MFCQ is satisfied.  Then, 
\[
\psi_{o}(\bar x)=\psi_{o}(\bar x;\bar H)=0
\]
if and only if $\bar x$ is a first order stationary point of \eqref{eq:nlp}; i.e. a KKT point for \eqref{eq:nlp}.  
\end{lem}
\bproof
First note that $\psi_{v}(x_k)=0$ is possible only if $d=0$ is a feasible solution for the quadratic subproblem \ref{eq:qpdircon}.  On the other hand, when $d=0$ is a feasible solution, $\psi_{o}(\bar x;\bar H)=0$ implies $d_{qp}=0$.  Finally, when $d_{qp}=0$, the KKT conditions for \eqref{eq:nlp} match the KKT conditions for \eqref{eq:qpdirobj}-\eqref{eq:qpdircon}.  Therefore, if $\psi_{v}(x_k)+\psi_{o}(\bar x;\bar H)=0$ implying $d_{qp}=0$, $c(\bar x)\leq 0$, and MFCQ is satisfied at $\bar x$, then $\bar x$ is a local solution for \eqref{eq:nlp}.  
Also, since the KKT conditions match, at a feasible regular (satisfying MFCQ) stationary point for \eqref{eq:nlp}, we have $d_{qp}=0$.  That implies $\psi_{o}(\bar x;\bar H)=0$ and $\psi_{v}(\bar x)=0$.

The same argument can be repeated for $\psi_{o}(\bar x)$.  
\eproof

\begin{lem}\label{lem:psiorel}
We have
\[
\psi_{o}(x_k) \leq \psi_{o}(x_k;H_k).
\] 
\end{lem}

\bproof 
This is a direct conclusion of the fact that
\[
\tilde q_f(d;x_k,H_k)\geq \tilde q_f(d;x_k,\sigma_n), \quad \forall d\in \tilde D;
\]
therefore,
\[
\min_{d\in \tilde D} \ \tilde q_f(d;x_k,H_k)\geq \min_{d\in \tilde D} \ \tilde q_f(d;x_k,\sigma_n).
\]
\eproof

Lemma~\ref{lem:psiorel} is important as it implies that a critical region result stated in terms of $\psi_{o}(x_k;H_k)$ provides the means to state a similar result in terms of $\psi_{o}(x_k)$.  We will next observe the effect of noise to the resulting optimality measure computed at a given iterate $x_k$.  To be able to do that, we will need to consider the regularity of the QP subproblems.

\begin{cond}[Regularity Condition]\label{cond:qpreg}
There exists a global bound $\Lambda_{qp} >0 $ for the optimal multipliers, that is for all $k$,
\[
\|\tilde\lambda_{k}\|_\infty\leq\Lambda_{qp}\quad \mbox{and} \quad \|\lambda_k\|_\infty\leq\Lambda_{qp}
\]
where $\tilde\lambda_{k}\in\R^m$ and $\lambda_{k}\in\R^m$
are the Lagrange multipliers for the solutions to the quadratic programming problems 
\begin{equation}\label{eq:qps}
\min_{d\in \tilde D} \ \tilde q_f(d;x_k,H_k)  \quad \mbox{and} \quad \min_{d\in D} \ q_f(d;x_k,H_k),
\end{equation}
respectively.
\end{cond}

The following is an implication of the well-known result on exact penalty functions (see for instance, \cite{nocedalwright},Theorem 17.3).
\begin{lem}\label{lem:qppen}
Suppose Condition~\ref{cond:qpreg} holds.  The unique solution $d_{qp}$ of the QP subproblem \eqref{eq:qpdirobj}-\eqref{eq:qpdircon} is a minimizer of the exact penalty function
\begin{equation}\label{eq:penqp}
    \tilde p_{qp}(d;x_k) = \tilde q_f(d;x_k,H_k) + m\Lambda_{qp} \|\max\{\tilde c_k+\tilde J_k^Td-r_v,0\}\|_\infty. 
\end{equation}
\end{lem}

\bproof
Define following sets regarding the solution $d_{qp}$.
\begin{align*}
\tilde{\mathcal{I}}_k = \{i: [\tilde c_k+\tilde J_k^Td_{qp}]_i<r_v\}, &\\
\tilde{\mathcal{A}}_k = \{i: [\tilde c_k+\tilde J_k^Td_{qp}]_i=r_v\}. &\\
\end{align*}
Since $d_{qp}$ is feasible for \eqref{eq:qpdircon}, we have $\|\max\{\tilde c_k+\tilde J_k^Td_{qp}-r_v,0\}\|_\infty=0$, and the subdifferential of $\|\max\{\tilde c_k+\tilde J_k^Td-r_v,0\}\|_\infty$ at $d_{qp}$ is given as the convex hull of the union of the subdifferentials of $\max\{[\tilde c_k+\tilde J_k^Td_{qp}-r_v]_i,0\}, i\in\mathcal{A}_k$.  Therefore, the subdifferential of $\tilde m_{qp}(d;x_k)$ at $d_{qp}$ can be written as 
\begin{equation*}
\partial \tilde p_{qp}(d_{qp}; x_k) = \tilde g_k + H_kd_{qp} + m\Lambda_{qp} \sum_{i\in\tilde{\mathcal{A}}_k}{w_i z_i[\tilde J_k]_{[i,:]}^T} 
\end{equation*}
with $w_i, z_i \in [0,1]$ and $\sum_i w_i = 1$.  

Consider optimal multipliers $\tilde \lambda_k$ of  \eqref{eq:qpdirobj}-\eqref{eq:qpdircon}.  So, for the choices $w_i=1/m, \forall i$ and $z_i=[\tilde \lambda_k]_i/\Lambda_{qp}$, the dual feasibility condition for \eqref{eq:qpdirobj}-\eqref{eq:qpdircon} at $d_{qp}$ implies that 
\[
0 \in \partial \tilde p_{qp}(d_{qp}; x_k),
\]
which in turn implies $d_{qp}$ is a minimizer of $\tilde p_{qp}(d; x_k)$.
\eproof

\bigskip

Before stating the next result, recall that 
\[
d_{qp} = \mbox{argmin}_{d\in \tilde D} \ \tilde q_f(d;x_k,H_k).
\]
For the non-noisy counterpart, we denote the corresponding unique solution by $\bar d_{qp}$; i.e.
\[
\bar d_{qp} = \mbox{argmin}_{d\in D} \ q_f(d;x_k,H_k).
\]
Following Lemma~\ref{lem:qppen}, it is easy to observe that $\bar d_{qp}$ is a minimizer of $p_{qp}(d;x_k)$ defined as
\begin{equation}\label{eq:penqpbar}
p_{qp}(d;x_k) = q_f(d;x_k,H_k) + m\Lambda_{qp} \|\max\{c_k+ J_k^Td-\bar r_v,0\}\|_\infty.
\end{equation}
The next lemma establishes a bound on the difference in between the minimum values of \eqref{eq:penqp} and \eqref{eq:penqpbar}.

\begin{lem}\label{lem:qpmoddiff}
Suppose Assumptions~\ref{assum:nlp} hold.  Then,
\[
|\tilde p_{qp}(d_{qp};x_k) - p_{qp}(\bar d_{qp};x_k)|\leq \epsilon_f + \epsilon_g r_l+m\Lambda_{qp}(2\epsilon_c+\epsilon_J (r_l+\Delta)).
\]
\end{lem}
\bproof
As we specified the formulation of \eqref{eq:lpfeas} using $\ell_\infty$ norm, we can write $r_v = \rho \textbf{1}$ and $\bar r_v = \bar \rho \textbf{1}$.
First we will observe
\begin{equation}\label{eq:diffvqp}
\left| \|\max\{c_k+J_k^Td_{qp}-\bar r_v,0\}\|_\infty - 
\|\max\{\tilde c_k+\tilde J_k^Td_{qp}-r_v,0\}\|_\infty \right| \leq 2\epsilon_c + \epsilon_Jr_l+\epsilon_J\Delta.
\end{equation}
Note that
\begin{align*}
&\|\max\{c_k+J_k^Td_{qp}-\bar r_v,0\}\|_\infty - 
\|\max\{\tilde c_k+\tilde J_k^Td_{qp}-r_v,0\}\|_\infty\\
& \leq \|\max\{c_k+J_k^Td_{qp}-\bar r_v,0\} - 
\max\{\tilde c_k+\tilde J_k^Td_{qp}-r_v,0\}\|_\infty \\
& \leq \|[c_k+J_k^Td_{qp}-\bar r_v] - 
[\tilde c_k+\tilde J_k^Td_{qp}-r_v]\|_\infty\\
& \leq \|c_k-\tilde c_k\|_\infty + \|J_k^Td_{qp}-\tilde J_k^Td_{qp}\|_\infty + |\bar \rho-\rho|\\
& \leq \|c_k-\tilde c_k\|_\infty + \|J_k^T-\tilde J_k^T\|_\infty\|d_{qp}\|_\infty + |\bar \rho-\rho|\\
& \leq \epsilon_c + \epsilon_J\|d_{qp}\| + |\bar \rho-\rho| \leq \epsilon_c + \epsilon_Jr_l + (\epsilon_c+\epsilon_J\Delta)
\end{align*}
by using Lemma \ref{lem:lendqp} and equation \eqref{eq:lvbound}.  Similarly we also have,
\[
\|\max\{\tilde c_k+\tilde J_k^Td_{qp}-r_v,0\}\|_\infty
-\|\max\{c_k+J_k^Td_{qp}-\bar r_v,0\}\|_\infty \leq \epsilon_c + \epsilon_Jr_l + (\epsilon_c+\epsilon_J\Delta)
\]
yielding \eqref{eq:diffvqp}.  Then, following similar lines to that of Lemma~\ref{lem:diffpsiv} we get
\begin{align*}
|p_{qp}(\bar d_{qp};x_k)-\tilde p_{qp}(d_{qp};x_k)| &= |(p_{qp}(\bar d_{qp};x_k)-p_{qp}(d_{qp};x_k))+(p_{qp}( d_{qp};x_k)-\tilde p_{qp}(d_{qp};x_k))|\\
& \leq |p_{qp}( d_{qp};x_k)-\tilde p_{qp}(d_{qp};x_k)|\\
& \leq | f_k + g_k^T d_{qp}+m\Lambda_{qp}\|c_k+J_k^Td_{qp}-r_v,0\}\|_\infty \\
& \qquad\qquad - \left[\tilde f_k + \tilde g_k^T d_{qp} + m\Lambda_{qp}\|\max\{\tilde c_k+\tilde J_k^Td_{qp}-\bar r_v,0\}\|_\infty\right] |\\
& \leq \epsilon_f + \epsilon_gr_l+m\Lambda_{qp}(2\epsilon_c + \epsilon_Jr_l+\epsilon_J\Delta).
\end{align*}
\eproof

The next lemma provides a bound on the gap between $\tilde\psi_o(x_k)$ and $\psi_o(x_k)$, similar to the result in Lemma~\ref{lem:diffpsiv} regarding $\tilde\psi_v(x_k)$ and $\psi_v(x_k)$.

\begin{lem}\label{lem:psiodiff}  Let Assumptions~\ref{assum:nlp} and Condition~\ref{cond:qpreg} hold.  Then,
\[
|\tilde\psi_o(x_k;H_k)-\psi_o(x_k;H_k)|\leq E_\epsilon^\psi
\]
for
\[
E_\epsilon^\psi = 2\epsilon_f+\epsilon_g r_l+m\Lambda_{qp}(2\epsilon_c+\epsilon_J (r_l+\Delta)).
\]
\end{lem}
\bproof
First note that since $d_{qp}$ is the optimal solution to the QP \eqref{eq:qpdirobj}-\eqref{eq:qpdircon}, 
\[
\|\max\{[\tilde c_k+\tilde J_k^Td_{qp}-r_v],0\}\|_\infty=0.
\]
Therefore, $\tilde p_{qp}(d_{qp};x_k)=\tilde q_f(d_{qp};x_k,H_k)$.  Similarly, we have $p_{qp}(\bar d_{qp};x_k)=q_f(\bar d_{qp};x_k,H_k)$.  Then, in the view of Lemma~\ref{lem:qpmoddiff},
\begin{align*}
|\tilde\psi_o(x_k;H_k)-\psi_o(x_k;H_k)| & = |(\tilde f_k-\tilde q_f(d_{qp};x_k,H_k)) - (f_k-q_f(\bar d_{qp};x_k,H_k))|\\ 
& = |(\tilde f_k-\tilde p_{qp}(d_{qp};x_k)) - (f_k-p_{qp}(\bar d_{qp};x_k))| \\
& \leq \epsilon_f + |\tilde p_{qp}(d_{qp};x_k) - p_{qp}(\bar d_{qp};x_k)|\\
& \leq \epsilon_f + \epsilon_f + \epsilon_g r_l+m\Lambda_{qp}(2\epsilon_c+\epsilon_J (r_l+\Delta))\\
& = 2\epsilon_f+\epsilon_g r_l+m\Lambda_{qp}(2\epsilon_c+\epsilon_J (r_l+\Delta)).
\end{align*}
\eproof

\paragraph{Line search. }

Before we proceed with the second (stronger) bound on the steplength, we establish the links between the merit function and the optimality measure.  The next result links the quadratic merit function model to the optimality measure. 
\begin{lem}\label{lem:psiomerit}
Suppose $\pi_k=\pi<\infty$ satisfies \eqref{piupdate} with $\theta_1\in(0,1)$. Then, the noisy quadratic model decrease of the merit function satisfies
\begin{equation}
\tilde q_\pi(0;x_k) - \tilde q_\pi(d_{qp};x_k) \geq \theta_1^2 [\tilde\psi_o(x_k;H_k)+\pi\tilde\psi_v(x_k)].
\end{equation}
\end{lem}

\bproof
\begin{align*}
\tilde q_\pi(0;x_k) - \tilde q_\pi(d_{qp};x_k) = -(\tilde g^Td_{qp} + \frac{1}{2}d_{qp}^TH_kd_{qp}) + \pi [\tilde l_v(0;x_k) - \tilde l_v(d_{qp};x_k)]
\end{align*}
If $\tilde g^Td_{qp} + \frac{1}{2}d_{qp}^TH_kd_{qp}<0$, then we immediately get
\begin{align*}
\tilde q_\pi(0;x_k) - \tilde q_\pi(d_{qp};x_k) &\geq |\tilde g^Td_{qp} + \frac{1}{2}d_{qp}^TH_kd_{qp}|+\pi [\tilde l_v(0;x_k) - \tilde l_v(d_{qp};x_k)]\\
& = \tilde \psi_o(x_k;H_k)+\pi\tilde\psi_v(x_k)\geq \theta_1^2 [\tilde\psi_o(x_k;H_k)+\pi\tilde\psi_v(x_k)]
\end{align*}
as $\theta_1\in(0,1)$.
On the other hand, if $\tilde g^Td_{qp} + \frac{1}{2}d_{qp}^TH_kd_{qp}>0$, we observe the fact that \eqref{piupdate} implies
\[
|\tilde g^Td_{qp} + \frac{1}{2}d_{qp}^TH_kd_{qp}| = \tilde g^Td_{qp} + \frac{1}{2}d_{qp}^TH_kd_{qp} \leq (1-\theta_1)\pi[\tilde l_v(0;x_k) - \tilde l_v(d_{qp};x_k)].
\] 
So,
\begin{align*}
\tilde q_\pi(0;x_k) - \tilde q_\pi(d_{qp};x_k) & = (\theta_1-(1+\theta_1))|\tilde g^Td_{qp} + \frac{1}{2}d_{qp}^TH_kd_{qp}|+\pi [\tilde l_v(0;x_k) - \tilde l_v(d_{qp};x_k)]\\
&\geq \theta_1|\tilde g^Td_{qp} + \frac{1}{2}d_{qp}^TH_kd_{qp}|+[-(1+\theta_1)(1-\theta_1)+1]\pi[\tilde l_v(0;x_k) - \tilde l_v(d_{qp};x_k)]\\
&=\theta_1|\tilde g^Td_{qp} + \frac{1}{2}d_{qp}^TH_kd_{qp}|+\theta_1^2\pi[\tilde l_v(0;x_k) - \tilde l_v(d_{qp};x_k)]\\
& \geq \theta_1^2|\tilde g^Td_{qp} + \frac{1}{2}d_{qp}^TH_kd_{qp}|+\theta_1^2\pi[\tilde l_v(0;x_k) - \tilde l_v(d_{qp};x_k)] \\
& = \theta_1^2 [\tilde\psi_o(x_k;H_k)+\pi\tilde\psi_v(x_k)]
\end{align*}
for $\theta_1\in(0,1)$.  This completes the proof.
\eproof

The following lemma is a key result for the analysis of line search, and for establishing a critical region based on the optimality measure.
\begin{lem}\label{lem:quadimp}
Let Assumptions~\ref{assum:nlp}, \ref{assum:alg}, and Condition~\ref{cond:qpreg} hold. 
Let $\beta \in (0,1)$. Then for any iterate $x_k$ such that
\begin{equation}\label{eq:psivcond}
\frac{\beta}{2} [\psi_o(x_k;H_k) + \pi\psi_v(x_k)] \geq \max\{E_\epsilon^\psi , \frac{1}{\theta_1^2}\epsilon_g r_l\}+\pi \max \{E_\epsilon^v , \frac{1}{\theta_1^2}(2\epsilon_c +\epsilon_Jr_l)\}
\end{equation}
we have 
\[
q_\pi(0;x_k) - q_\pi(d_{qp};x_k) \geq \theta_3 [\psi_o(x_k;H_k)+\pi\psi_v(x_k)],
\]
where $\theta_3= (1-\beta)\theta_1^2$.
If we choose $\beta =\frac{1}{2}$, then
$\theta_3=\frac{1}{2}\theta_1^2$.
\end{lem}

\bproof
By Lemma~\ref{lem:deltabound},
\begin{align}
q_\pi(0;x_k) - q_\pi(d_{qp};x_k) &= \tilde q_\pi(0;x_k) - \tilde q_\pi(d_{qp};x_k) + \delta_k \nonumber\\
&\geq \tilde q_\pi(0;x_k) - \tilde q_\pi(d_{qp};x_k)-\epsilon_g r_l - 2\pi\epsilon_c -\pi\epsilon_Jr_l. \label{ddelta}
\end{align}
Also, using Lemmas~\ref{lem:diffpsiv} and \ref{lem:psiodiff} we get,
\[
\tilde\psi_o(x_k;H_k)+\pi\tilde\psi_v(x_k) \geq \psi_o(x_k;H_k)+\pi\psi_v(x_k) - E_\epsilon^{\psi}-\pi E_\epsilon^v. 
\]
Therefore, in the view of Lemma~\ref{lem:psiomerit}, \eqref{ddelta} yields
\begin{align*}
q_\pi(0;x_k) - q_\pi(d_{qp};x_k) & \geq \theta_1^2[\tilde\psi_o(x_k;H_k)+\pi\tilde\psi_v(x_k)] -\epsilon_g r_l - 2\pi\epsilon_c -\pi\epsilon_Jr_l\\
& \geq \theta_1^2[\psi_o(x_k;H_k)+\pi\psi_v(x_k)]-(\theta_1^2E_\epsilon^{\psi}+\epsilon_gr_l)- \pi (\theta_1^2E_\epsilon^v+2\epsilon_c+\epsilon_Jr_l)\\
& \geq \theta_1^2[\psi_o(x_k;H_k)+\pi\psi_v(x_k)]-\theta_1^2 \beta \psi_o(x_k;H_k)-\pi \theta_1^2 \beta \psi_v(x_k)\\
&= (1-\beta)\theta_1^2[\psi_o(x_k;H_k)+\pi\psi_v(x_k)].
\end{align*}
\eproof

Here, the particular choice of $\beta=\frac{1}{2}$ is arbitrary and given to simplify the analysis; all derivations would work for any $\beta\in(0,1)$.  Next, we will show that for all iterates satisfying \eqref{eq:psivcond} there is a positive lower bound for $\alpha_k$.  

\begin{lem}\label{lem:alpha}
Define $L_\pi$ as in Lemma~\ref{lem:alpha_weak}.  If \eqref{eq:psivcond} holds for $\beta=\frac{1}{2}$, then the line search condition \eqref{eq:lscond} is satisfied for any positive $\alpha$ satisfying 
\begin{equation}\label{eq:alpha}
 \alpha\leq \frac{2(1-\theta_4)(1-\theta_2)\theta_3 [\psi_o(x_k;H_k)+\pi\psi_v(x_k)]}{L_\pi \|d_{qp}\|^2}+\frac{\sigma_n}{L_\pi}.
\end{equation}
for $\theta_3>0$ specified in Lemma~\ref{lem:quadimp} and for 
\[
\theta_4 =\frac{1}{2}\frac{\theta_2}{1-\theta_2}.
\]
\end{lem}

\bproof
Since \eqref{eq:psivcond} implies 
\[
-\delta\leq\epsilon_g r_l + \pi(2\epsilon_c +\epsilon_Jr_l) \leq \frac{\theta_1^2}{4} [\psi_o(x_k;H_k) + \pi\psi_v(x_k)],
\]
inserting the definitions of $\theta_3 = \frac{\theta_1^2}{2}$ and $\theta_4=\frac{1}{2}\frac{\theta_2}{1-\theta_2}$ yields
\[
-\delta \leq \frac{\theta_3}{2} [\psi_o(x_k;H_k) + \pi\psi_v(x_k)] = \theta_3\theta_4\frac{(1-\theta_2)}{\theta_2}[\psi_o(x_k;H_k) + \pi\psi_v(x_k)].
\]
So, we get
\begin{equation}\label{eq:delta}
\theta_2\delta \geq -\theta_4(1-\theta_2)\theta_3[\psi_o(x_k;H_k)+\pi\psi_v(x_k)].
\end{equation}  

\bigskip
Following similar lines to the proof of Lemma~\ref{lem:alpha_weak} we have 
\[
\phi_{\pi}(x_k) - \phi_{\pi}(x_k+\alpha d_{qp}) \geq \alpha[q_\pi(0;x_k) - q_\pi(d_{qp};x_k)]+\frac{\alpha}{2}d_{qp}H_kd_{qp}-\frac{L_\pi}{2}\alpha^2 \|d_{qp}\|^2.
\]
By Lemma~\ref{lem:quadimp}, this yields
\begin{align*}
\phi_{\pi}(x_k) - \phi_{\pi}(x_k+\alpha d_{qp}) 
& \geq \alpha[\theta_2+(1-\theta_2)][q_\pi(0;x_k) - q_\pi(d_{qp};x_k)]+\alpha\frac{\sigma_n}{2}\|d_{qp}\|^2-\frac{L_\pi}{2}\alpha^2 \|d_{qp}\|^2\\
& \geq \alpha\theta_2[\tilde q_\pi(0;x_k) - \tilde q_\pi(d_{qp};x_k)+\delta]+\alpha(1-\theta_2)\theta_3 [\psi_o(x_k;H_k)+\pi\psi_v(x_k)]\\
& \qquad +\alpha\frac{\sigma_n}{2}\|d_{qp}\|^2-\frac{L_\pi}{2}\alpha^2\|d_{qp}\|^2
\end{align*}
implying 
\begin{align*}
\tilde\phi_{\pi}(x_k) - \tilde\phi_{\pi}(x_k+\alpha d_{qp}) & \geq 
\alpha\theta_2[\tilde q_\pi(0;x_k) - \tilde q_\pi(d_{qp};x_k)] -2\epsilon_R +\alpha\theta_2\delta \\
& \ \ +\alpha(1-\theta_2)\theta_3 [\psi_o(x_k;H_k)+\pi\psi_v(x_k)] 
 +\alpha\frac{\sigma_n}{2}\|d_{qp}\|^2-\frac{L_\pi}{2}\alpha^2\|d_{qp}\|^2
\end{align*}
by \eqref{eq:eqrelax}.  Therefore, the line search condition \eqref{eq:lscond} holds for $\alpha$ satisfying
\begin{equation}\label{eq:lsguarantee}
\alpha\theta_2\delta+\alpha(1-\theta_2)\theta_3 [\psi_o(x_k;H_k)+\pi\psi_v(x_k)]+\alpha\frac{\sigma_n}{2}\|d_{qp}\|^2-\frac{L_\pi}{2}\alpha^2\|d_{qp}\|^2 \geq 0.
\end{equation}
By ~\eqref{eq:delta}, the inequality \eqref{eq:lsguarantee} is guaranteed to hold when 
\[
h(\alpha) := (1-\theta_4)(1-\theta_2)\theta_3 [\psi_o(x_k;H_k)+\pi\psi_v(x_k)]\alpha+\alpha\frac{\sigma_n}{2}\|d_{qp}\|^2-\frac{L_\pi}{2}\alpha^2 \|d_{qp}\|^2 \geq 0
\] 
holds.  The concave quadratic function $h(\alpha)$ takes nonnegative values in between its two roots; $\alpha=0$ and
\[
\alpha = \frac{2(1-\theta_4)(1-\theta_2)\theta_3 [\psi_o(x_k;H_k)+\pi\psi_v(x_k)] +\sigma_n\|d_{qp}\|^2}{L_\pi \|d_{qp}\|^2} \equiv \bar{\alpha}_k,
\]
which always exists.  
\eproof

Now we can use this result to establish a positive lower bound on $\alpha_k$. It is worth noting that, although the bound, $\bar\alpha $,
depends on the noise bounds of Assumptions \ref{assum:nlp}, it is bounded below by $\frac{\sigma_n}{2L_\pi}$. 
\begin{lem}\label{lem:baralpha}
At any iterate $x_k$ where \eqref{eq:psivcond} holds the steplength $\alpha_k$ provided by Algorithm~\ref{alg:rsqp} satisfies 
\[
\alpha_k\geq \bar\alpha.
\]
where
\begin{equation}        \label{eq:alphabar}
\bar\alpha=\frac{4(1-\theta_4)(1-\theta_2)\theta_3 }{ r_{\ell}^2} \min \{ \frac{\max\{E_\epsilon^\psi ,\epsilon_g r_l / \theta_1^2 \}}{L_f} , \frac{\max\{E_\epsilon^v , (2\epsilon_c +\epsilon_Jr_l)/\theta_1^2\} }{L_c}\}+\frac{\sigma_n}{2L_\pi}.
\end{equation}
\end{lem}
\bproof
Note that equation \eqref{eq:alpha} is implied by the stronger condition
\[
\alpha \leq  \frac{2(1-\theta_4)(1-\theta_2)\theta_3 [\psi_o(x_k;H_k)+\pi\psi_v(x_k)]}{( L_f+\pi L_c) r_{\ell}^2}+\frac{\sigma_n}{L_\pi} , 
\]
using the bound in Lemma \ref{lem:lendqp}, and
since $L_\pi = L_f+\pi L_c$.

Note also that by \eqref{eq:psivcond} this equation is implied by the stronger condition
\begin{equation}  \label{pifunc}
\alpha \leq  \frac{8(1-\theta_4)(1-\theta_2)\theta_3 }{( L_f+\pi L_c) r_{\ell}^2} [\max\{E_\epsilon^\psi , \frac{1}{\theta_1^2}\epsilon_g r_l\}+\pi \max \{E_\epsilon^v , \frac{1}{\theta_1^2}(2\epsilon_c +\epsilon_Jr_l)\}]+\frac{\sigma_n}{L_\pi}.
\end{equation}
It is straightforward to see that the right hand side of \eqref{pifunc} is either a monotone increasing or monotone decreasing function of $\pi$ on the interval $[0,\infty)$. It follows that this quantity is bounded below by its value at $\pi=0$ or by its limiting value as $\pi \rightarrow \infty$.

Therefore it follows from Lemma \ref{lem:alpha} that for any positive value of $\pi$,  if \eqref{eq:psivcond} holds at $x_k$, then line search condition \eqref{eq:lscond} is satisfied for all $\alpha$ satisfying
\[
\alpha \leq  \frac{8(1-\theta_4)(1-\theta_2)\theta_3 }{ r_{\ell}^2} \min \{ \frac{\max\{E_\epsilon^\psi ,\epsilon_g r_l / \theta_1^2 \}}{L_f} , \frac{\max\{E_\epsilon^v , (2\epsilon_c +\epsilon_Jr_l)/\theta_1^2\} }{L_c} \}+\frac{\sigma_n}{L_\pi}.
\]
Now if the backtracking line search  in Step 5 of Algorithm \ref{alg:rsqp} chooses $\alpha_k< 1$ it must be because \eqref{eq:lscond} was violated at $\alpha=2\alpha_k$, which by the above would imply that $2\alpha_k$ is greater than the right hand side of \eqref{pifunc}.  It follows that for all $x_k$ satisfying \eqref{eq:psivcond}, $\alpha_k \geq \bar\alpha$, where $\bar\alpha$ is given by \eqref{eq:alphabar}.

\eproof

\paragraph{Penalty parameter.}

\begin{lem}\label{lem:penalty}
Suppose the sequence produced by Algorithm 1 satisfies Condition~\ref{cond:qpreg}. Then for a sufficiently large finite value of the penalty parameter, condition \eqref{piupdate} is satisfied for any $k$.  It follows that there exists $\hat{\pi}<\infty$ and an integer $k_1$ such that $\pi_k = \hat\pi$, for all $k \geq k_1$.
\end{lem}

\bproof
By Lemma~\ref{lem:qppen}, $\tilde p_{qp}(d_{qp};x_k)\leq \tilde p_{qp}(0;x_k)$.  Therefore,
\begin{align*}
0 & \geq \tilde g_k^Td_{qp} + \frac{1}{2}d_{qp}^TH_kd_{qp} + m\Lambda_{qp} \|_\infty\max\{[\tilde c_k+\tilde J_k^Td_{qp}-r_v],0\}\| - m\Lambda_{qp} \|\max\{[\tilde c_k-r_v],0\}\|_\infty\\
& = \tilde g_k^Td_{qp} + \frac{1}{2}d_{qp}^TH_kd_{qp} - m\Lambda_{qp} \|\max\{[\tilde c_k-r_v],0\}\|_\infty
\end{align*}
as $\|\max\{[\tilde c_k+\tilde J_k^Td_{qp}-r_v],0\}\|_\infty=0$ by optimality of $d_{qp}$ for \eqref{eq:qpdirobj}-\eqref{eq:qpdircon}.
Since \eqref{eq:lpfeas} is formulated using $\ell_\infty$ norm, we can write $r_v = \rho \textbf{1}$.  Then, 
\begin{align*}
\|\max\{[\tilde c_k-r_v],0\}\|_\infty & = \|\max\{\tilde c_k- \rho \textbf{1},0\}\|_\infty\\
& = \|\max\{\tilde c_k,0\}\|_\infty-\rho\\
& = \tilde l_v(0;x_k) - \tilde l_v(d_{qp};x_k).
\end{align*}
Here, $\|\max\{\tilde c_k- \rho \textbf{1},0\}\|_\infty=\|\max\{\tilde c_k,0\}\|_\infty-\rho$ holds because 
\[
\rho = \min_{\|d\|\leq\Delta} \|\max\{\tilde c_k+\tilde J_kd,0\}\|_\infty \leq \|\max\{\tilde c_k,0\}\|_\infty.
\]
Gathering all we get
\begin{align*}
\tilde g_k^Td_{qp} + \frac{1}{2}d_{qp}^TH_kd_{qp} & \leq m\Lambda_{qp} \|\max\{[\tilde c_k-r_v],0\}\|_\infty\\
&\leq m\Lambda_{qp} (\tilde l_v(0;x_k) - \tilde l_v(d_{qp};x_k)).
\end{align*}
This implies
\begin{align*}
\tilde q_{\pi}(0;x_k) - \tilde q_{\pi}(d_{qp};x_k)] \geq \theta_1 \pi[\tilde l_v(0;x_k) - \tilde l_v(d_{qp};x_k)]
\end{align*}
is satisfied for any $\pi \geq \frac{1}{1-\theta_1}m\Lambda_{qp}$,
which implies that $\{\pi_k \}$ is bounded. 
Then since any change in $\pi_k$ in Step 4 of the algorithm increases it by at least $0.1 \pi_k$, eventually $\pi_k$ will stop increasing, and remain fixed at a value $\hat{\pi} \leq \frac{1}{1-\theta_1}m\Lambda_{qp}$. 
\eproof

\paragraph{A stronger feasibility result.} Thanks to the stronger steplength bound obtained in Lemma~\ref{lem:baralpha}, it is possible to give a stronger result for feasibility when regularity holds.

To simplify notation, let us define
\[
E_\epsilon(\pi)=\frac{1}{\pi}\max\{E_\epsilon^\psi , \frac{1}{\theta_1^2}\epsilon_g r_l\}+ \max \{E_\epsilon^v , \frac{1}{\theta_1^2}(2\epsilon_c +\epsilon_Jr_l)\}
\]
so that \eqref{eq:psivcond} with $\beta=\frac{1}{2}$ can be represented as
\[
\frac{1}{4} [\frac{1}{\pi}\psi_o(x_k;H_k) + \psi_v(x_k)]\geq E_\epsilon(\pi).
\]

\begin{lem}\label{cor:rfeas}
Suppose Assumptions~\ref{assum:nlp} and \ref{assum:alg} hold. Suppose also that Condition~\ref{cond:qpreg} holds for the iterates of Algorithm~\ref{alg:rsqp}.  Then, the algorithm has infinitely many iterates in a region $\mathcal{F}_o$ defined as 
\begin{equation}\label{eq:crit3}
\mathcal{F}_o = \left\{x: \psi_v(x) \leq \max\left( 4 E_\epsilon(\hat \pi), \frac{32}{3\pi\theta_1\theta_2\bar\alpha}\epsilon_R\right)\right\}.
\end{equation}
\end{lem}
\bproof
By Lemma~\ref{lem:penalty}, $\pi_k = \hat \pi$ for $k$ sufficiently large.  Then, for $x_k\notin \mathcal{F}_o$ it holds that
\[
\psi_v(x_k)\geq  4E_\epsilon(\hat \pi)\geq 4E_\epsilon^v \quad \Rightarrow \quad \psi_v(x_k)-E_\epsilon^v \geq  \frac{3}{4}\psi_v(x_k),
\]
and
\[
\psi_v(x_k)\geq \frac{4}{3}\frac{8}{\hat \pi\theta_1\theta_2\bar\alpha}\epsilon_R \quad \Rightarrow \quad \frac{4}{\hat \pi}\epsilon_R \leq \theta_1\theta_2\bar\alpha \frac{3}{8}\psi_v(x_k).
\]
Recall the scaled merit function $\phi_\pi^s$ defined in \eqref{eq:scaledmerit}.  Following the same lines as in the proof of Theorem~\ref{thm:crfeas}, and using the lower bound on $\alpha_k$ in Lemma~\ref{lem:baralpha}, we obtain 
\begin{align*}
\phi^s_{\hat \pi}(x_k) - \phi^s_{\hat\pi}(x_{k+1}) & \geq \theta_1\theta_2\bar\alpha(\psi_v(x_k)-E_\epsilon^v) - \frac{4}{\hat \pi}\epsilon_R\\
& \geq \frac{3}{8}\theta_1\theta_2\bar\alpha \psi_v(x_k)
\end{align*}
for $x_k\notin \mathcal{F}_o$.  If $\{x_k\}$ has only finite number of elements in $\mathcal{F}_o$, then $\psi_v(x_k)$ cannot stay bounded away from zero because $\phi^s_{\pi}$ is bounded from below.
\eproof
The result provided in Corollary~\ref{cor:rfeas} is stronger than that of Theorem~\ref{thm:crfeas} as the bound on $\psi_v$ is proportional to the noise level rather than its square-root.

\paragraph{Convergence to a critical region for optimality.} We next present a result providing a threshold on $\psi_o(x_k)+\pi\psi_v(x_k)$ as a function of the noise level, above which the steps of Algorithm~\ref{alg:rsqp} guarantee to improve the (non-noisy) merit function.  That gives a critical region for optimality.

\begin{thm}\label{thm:cropt} 
Suppose all assumptions of Lemma~\ref{cor:rfeas} hold. Then, 
the algorithm has infinitely many iterates in a region $\mathcal{C} \subseteq{\mathcal F_o}$ (a critical region for optimality) specified as 
\begin{equation}\label{eq:crit2}
\mathcal{C} = \left\{x: \psi_o(x)+\hat\pi\psi_v(x) \leq \max\left( 4\hat\pi E_\epsilon(\hat\pi), \frac{32}{3\theta_1\theta_2\bar\alpha}\epsilon_R\right)\right\}.
\end{equation}
Each time the iterates leave $\mathcal{C}$, they will return to $\mathcal{C}$.
\end{thm}

\bproof
Similar to the proofs of Theorem~\ref{thm:crfeas} and Lemma~\ref{cor:rfeas}, for $x_k\notin C$ we get
\begin{align*}
 \tilde\phi_{\hat \pi}(x_k) - \tilde\phi_{\hat\pi}(x_{k+1}) &\geq \theta_2\alpha_k[\tilde q_{\hat\pi}(0;x_k) - \tilde q_{\hat\pi}(d_{qp};x_k)] - 2\epsilon_R \\
\Rightarrow \ \phi_{\hat\pi}(x_k) - \phi_{\hat\pi}(x_{k+1}) &\geq \theta_2\bar\alpha\theta_1^2[\psi_o(x_k;H_k)+\hat\pi \psi_v(x_k)-(E_\epsilon^\psi+\hat\pi E_\epsilon^v)] - 4\epsilon_R\\
  &\geq \theta_2\bar\alpha\theta_1^2[\psi_o(x_k;H_k)+\hat\pi \psi_v(x_k)-\hat\pi E_\epsilon(\hat\pi))] - 4\epsilon_R\\
  &\geq \theta_2\bar\alpha\theta_1^2\frac{3}{8}[\psi_o(x_k;H_k)+\hat\pi \psi_v(x_k)] \\
\Rightarrow \ \phi_{\hat\pi}(x_k) - \phi_{\hat\pi}(x_{k+1}) & \geq \theta_2\bar\alpha\theta_1^2\frac{3}{8}[\psi_o(x_k)+\hat\pi \psi_v(x_k)]
\end{align*}
by Lemmas~\ref{lem:psiomerit}, \ref{lem:psiodiff}, and \ref{lem:diffpsiv}.  Since $\phi_{\hat\pi}(x)$ is bounded from below, that implies $\psi_o(x_k)+\hat\pi \psi_v(x_k)$ cannot stay bounded away from zero for $x_k\notin C$.
\eproof

Theorem~\ref{thm:cropt} also provides an upper bound for the merit function value as $k\rightarrow \infty$.  Let 
\[
\phi_{\max}^{\mathcal C} = \max_{x_k\in \mathcal{C}} \phi_{\hat\pi}(x_k), \qquad \mbox{and} \qquad \bar \phi_{\max}= \phi_{\max}^{\mathcal C} + 2\epsilon_f+2\hat\pi\epsilon_c,
\]
where $\mathcal C$ is the set given by \eqref{eq:crit2}.  When the iterates leave $\mathcal C$ at any iteration $k$, the maximum value the merit function can take is 
\begin{align*}
\phi_{\max}^{\mathcal C}+2\epsilon_R 
= & \ \phi_{\max}^{\mathcal C}+2\epsilon_f+2\pi_k\epsilon_c\\
\leq & \ \phi_{\max}^{\mathcal C}+2\epsilon_f+2\hat\pi\epsilon_c = \bar \phi_{\max}.
\end{align*}
Therefore, we have $\underset{k\rightarrow \infty}{\limsup} \ \phi_{\hat\pi}(x_k) = \bar \phi_{\max}$, implying that the iterates $x_k$ of Algorithm~\ref{alg:rsqp} cannot leave the level set $\{x: \phi_{\hat\pi}(x)\leq \bar \phi_{\max}\}$ for all $k$ sufficiently large.
 
\subsection{Discussion of Regularity}

The findings of our work regarding optimality rely on Condition~\ref{cond:qpreg}.  The quantity $\psi_o(x)$ is a proper measure only when that regularity condition holds, and a critical region for optimality is well-defined only in this case.  

The regularity condition \eqref{cond:qpreg} is stated in terms of the multipliers of the quadratic programming subproblems \eqref{eq:qps}.  As these problems are feasible and convex by construction, boundedness of their multipliers is guaranteed if the interior of their feasible regions are nonempty -- i.e. if \emph{Slater's condition} holds.  On the other hand, it is not hard to see that Slater's condition for the QP feasible region \eqref{eq:D} constructed at $w$ is equivalent to the MFCQ condition for the NLP \eqref{eq:nlp} at $w$, 
\[
\exists d : \nabla c_i(w)^Td < 0 , i\in \{i: c_i(w)=0\},
\]
when $w$ is feasible for \eqref{eq:nlp}; that is, when $v(w)=r_v=0$.  So, when $v(w)=0$, the failure of Slater's condition for the set $D=\{d : c(w)+\nabla c(w)^Td\leq 0\}$ is equivalent to the failure of MFCQ for \eqref{eq:nlp} at $w$. 

Note also that Condition \ref{cond:qpreg} is closely related to the regularity condition (and the stability results) given in \cite{fiacco}.  Indeed, the first part of the regularity condition in \cite{fiacco} (see Definition 3.2) is always satisfied by the subproblems in \eqref{eq:qps} as they are always feasible and convex by construction; also, they are always bounded thanks to their strictly convex quadratic objectives.  On the other hand, the second part of the regularity condition of \cite{fiacco} regarding the boundedness of multipliers cannot fail in a sufficiently small neighborhood of a regular problem (see Theorem 3.5 of \cite{fiacco}).  In other words, the failure of Condition \ref{cond:qpreg} could be possible if any of the subproblems \eqref{eq:qps} is \emph{close to} violating Slater's condition.  But since the distance in between the constraint parameters of the two subproblems is bounded by $\epsilon_c$ and $\epsilon_J$, this is only possible if Slater's condition can fail via a perturbation of $O(\epsilon_c,\epsilon_J)$.  We elaborate on this in Definition~\ref{def:almostsingular} and Corollary~\ref{cor:regres}.  

\begin{Def}\label{def:almostsingular}
Problem \eqref{eq:nlp} is $\epsilon$-close to violation of MFCQ at a solution point $w$ satisfying $v(w)\leq\epsilon$ if there exists $\bar c\in\mathbb{R}^m, \bar J\in\mathbb{R}^{m\times n}$ such that
\[
\|\bar c-c(w)\|\leq \epsilon, \ \ \|\bar J-J(w)\|\leq \epsilon,
\]
and $\nexists d\in\mathbb{R}^n$ satisfying
\[
\bar c +\bar J^T d < r_v , \qquad r_v\in \underset{r\geq \bar c+ \bar J^Tv, \ \|v\|\leq \Delta}{\mbox{arg}\min} \|r\|.
\]
\end{Def}

Definition~\ref{def:almostsingular} expresses the existence of a perturbation of the feasible region $D$ defined in \eqref{eq:D} for which Slater's condition fails.  The following example illustrates $\epsilon$-closeness to violation of regularity.

\paragraph{Example.}  
The optimization problem instance is
\begin{align*}
\min_{x_1,x_2} \ \ & x_1 + x_2\\
& x_2^2 \geq a_2\\
& 0.5x_2^2 + x_1x_2 \geq 0\\
& x_1, x_2 \geq 0,
\end{align*}
with $a_2=10^{-4}$.  Without noise, rSQP correctly converges to the solution $x^\ast=(0.0,\sqrt{a_2})$, where the two active constraints have linearly independent gradients.  During the run of the algorithm, the value of $\pi$ is never increased.

We rerun the algorithm on this problem starting from the same initial solution, but this time injecting noise to function and gradient evaluations with $\epsilon_f=\epsilon_g=10^{-2}$.  Since the Jacobian entries have norm smaller than the noise level when the iterates are close to the above mentioned solution, we pick noise values such that the Jacobian becomes singular.  Eventually, the penalty parameter keeps growing.  In other words, even if MFCQ holds at $x^\ast$ for the non-noisy problem, the regularity condition \eqref{cond:qpreg} fails to hold in the presence of noise.  Indeed, in view of Definition~\ref{def:almostsingular}, this problem is $\epsilon$-close to violation of regularity at $x^\ast$ for $\epsilon=10^{-2}$.  That is, at $w=x^\ast$, 
\[
c(w)=\begin{bmatrix}
    0.0\\ 0.0\\ 0.0\\ -1e-2
\end{bmatrix} 
\quad \mbox{ and }\quad 
J(w)^T=\begin{bmatrix} 
0.0 & -2e-2\\ 
-1e-2 & -1e-2\\ 
-1 & 0\\ 
0 & -1
\end{bmatrix}. 
\]
It is possible to have $\tilde c_2(w)=1e-2$ and $\tilde \nabla c_2(w)^T=[0 \ \ 0]$.  $[r_v]_2=1e-2$ provides feasibility, but the strict inequality $\tilde c_2(w) + \tilde \nabla c_2(w)^Td < [r_v]_2$ cannot be satisfied for any $d$.

\bigskip

We next provide a result that explains convergence of the iterates of Algorithm~\ref{alg:rsqp} to points $\epsilon$-close to violation of regularity.

\begin{cor}\label{cor:regres}
Suppose Assumptions 2.1 and 3.1 hold and consider the iterates $\{x_k\}$ of Algorithm 1.  Assume further that the level sets of $\psi_v$ are compact. Suppose that the algorithm approaches feasibility in the sense that 
$$
\liminf v(x_k) \leq \epsilon.
$$
Then, if Condition~\ref{cond:qpreg} is not satisfied, there exists a subsequence of $\{x_k\}$ with an accumulation point $w$ such that \eqref{eq:nlp} is $\epsilon$-close to violation of MFCQ at $w$ in the sense of Definition~\ref{def:almostsingular} if $v(w)\leq \epsilon$.
\end{cor}

\bproof
As Condition~\ref{cond:qpreg} does not hold, there is a subsequence of iterations such that
\begin{equation}   \label{lambdaunbdd} 
\max\{\|\lambda_{k_j}\|, \|\tilde{\lambda}_{k_j}\|\} \rightarrow \infty .
\end{equation} 
This implies that at least either $\|\tilde{\lambda}_{k_j}\| \rightarrow \infty$ or $\|\lambda_{k_j}\| \rightarrow \infty$. We consider first the case where $\tilde{\lambda}_{k_j}$ is unbounded.

By the assumption of boundedness of the iterates, this subsequence must have a convergent subsequence such that $x_{k_j} \rightarrow w$ for some $w$. Additionally, since the noisy values and gradients at these iterates satisfy the noise bounds \eqref{eq:noisebounds1} and \eqref{eq:noisebounds2}, the sequences $\tilde{f}_{k_j}, \tilde{g}_{k_j},\tilde{c}_{k_j},$ and $\tilde{J}_{k_j}$ are bounded and must have a cluster point $\tilde{f}_{*}, \tilde{g}_{*},\tilde{c}_{*},$ and $\tilde{J}_{*}$ such that
\begin{equation}\label{eq:limitbounds1}
|\tilde f_* - f(w)|\leq \epsilon_f, \qquad \|\tilde g_*-\nabla f(w)\|\leq \epsilon_g,  
\end{equation}
and
\begin{equation}\label{eqlimitbounds2}
\|\tilde c_*-c(w)\|\leq \epsilon_c \qquad \| \tilde J_*-\nabla c(w)\| \leq \epsilon_J. 
\end{equation}

Now consider the QP defined by $\tilde{f}_{*}, \tilde{g}_{*},\tilde{c}_{*},$ and $\tilde{J}_{*}$,
\begin{align}\label{eq:qpdircon2}
\mbox{minimize}_d \ \ & \tilde f_* + \tilde g_*^Td + \frac{1}{2}d^TH_*d\\
s.t. \quad & \tilde c_*+\tilde J_*^Td \leq r_v.\nonumber
\end{align}
The QP \eqref{eq:qpdircon2} is feasible by construction.  If there were a vector $\bar{d}$ strictly feasible for \eqref{eq:qpdircon2} such that $\tilde c_*+\tilde J_*^T\bar d < r_v$, then the QP \eqref{eq:qpdircon2} would be regular according to Definition 3.2 of \cite{fiacco} (i.e. the Slater condition holds for \eqref{eq:qpdircon2}), which is true if and only if the multipliers at the solution of \eqref{eq:qpdircon2} are bounded \cite{gauvin1977necessary}.  Then, according to Theorem 3.5 of \cite{fiacco}, the same regularity condition would hold for any QP in a neighborhood of \eqref{eq:qpdircon2}, and all QPs obtained by perturbing \eqref{eq:qpdircon2} by a sufficiently small amount would have bounded multipliers.

This contradicts our assumption that the multipliers on the subsequence diverge.   Therefore no such $\bar{d}$ exists. It follows by Definition 3.22 that \eqref{eq:nlp} is $\epsilon$-close to violation of regularity at $w$. 

If we consider the second possibility for \eqref{lambdaunbdd}, an essentially identical argument leads to the same contradiction, thus establishing the result.
\eproof

\bigskip

Before closing this section, we should also mention the relevance of the results of our paper to that of \cite{gabo}.  When there is no noise; i.e. when $\epsilon_f=\epsilon_g=\epsilon_c=\epsilon_J=0$, it is not hard to see that the line of our analysis lead to results that coincide with that of Theorem 3.14 of \cite{gabo}.

\begin{cor}
Suppose $\epsilon_f=\epsilon_g=\epsilon_c=\epsilon_J=0$.  Then, 
\begin{itemize}
\item[a.] if Condition~\ref{cond:qpreg} holds, then accumulation points of the sequence $\{x_k\}$ produced by Algorithm~\ref{alg:rsqp} are either infeasible stationary points, or KKT points of \eqref{eq:nlp};  
\item[b.] otherwise if Condition~\ref{cond:qpreg} fails, there is a limit point $x^\ast$ that is an infeasible stationary point or a feasible limit point $x^\ast$ where MFCQ fails.
\end{itemize}
\end{cor}

\bproof
By Theorem~\ref{thm:crfeas}, we have $\psi_v(x^\ast)=0$ at any accumulation point $x^\ast$ when $\epsilon_g=\epsilon_c=\epsilon_J=0$.  Then, if $x^\ast$ is infeasible, it is an infeasible stationary point.  If Condition~\ref{cond:qpreg} holds and $x^\ast$ is feasible, then by Theorem~\ref{thm:cropt} we have $\psi_o(x^\ast)=0$ implying that $x^\ast$ is a KKT point.  Otherwise, if Condition~\ref{cond:qpreg} does not hold, there are two cases: either $\{x_k\}$ has an infeasible limit point, or all limits points of $\{x_k\}$ are feasible.  But by Corollary~\ref{cor:regres}, when Condition~\ref{cond:qpreg} does not hold $\{x_k\}$ has a subsequence with a limit point $x^\ast$ such that if $x^\ast$ is feasible then it is a limit point where MFCQ fails.   
\eproof

\section{Numerical Experiments}

In this section, we conduct preliminary numerical experiments to observe the average behavior of Algorithm~\ref{alg:rsqp} in getting close to stationary points of feasibility / optimality measures.  We conduct experiments on a subset of small scale CUTEr instances for this purpose.  We first solve the instances without noise to get reference solution values.  Then we solve the same instances by injecting random noise at various levels, and observe the effect of noise on the resulting solution set.

We should point out that the non-noisy version of the rSQP algorithm itself is new, although it has close similarities to existing algorithms.  Therefore, there is no extensive experimental results available on the performance of this algorithm in the literature.  However, such an extensive testing is out of the scope of this paper.  Here we concentrate only on robustness in the presence of noise.

We implemented the algorithm in Python by employing the \texttt{scipy.optimize} and \texttt{qpsolvers} packages for solving the subproblems.  For practical convenience, we solve subproblem \eqref{eq:lpfeas} as an LP by restricting the length of $d$ via its $\ell_\infty$ norm rather than $\ell_2$. In our implementation, we employ BFGS updates to get $H_k$, which ensures positive definiteness.  Since $H_k$ is ideally an approximation to the Hessian of the Lagrangian, we use curvature pairs $s_k=\alpha_kd_{qp}$ and $y_k=(\tilde g_{k+1}+\lambda_{k+1} \tilde J_{k+1})-(\tilde g_{k}+\lambda_{k+1} \tilde J_k)$.  Here $\lambda_{k+1}$ is computed as the least-squares multipliers.  The BFGS updates are done only if the following two conditions
\[
y_k^Ts_k \geq 10^{-3}s_k^TH_ks_k \qquad \mbox{and} \qquad y_k^Ts_k \geq (\epsilon_g +\|\lambda_k\|_\infty\epsilon_J)\|s_k\|
\]
hold.  Otherwise, the updates are skipped.  The second rule above ensures the curvature condition is satisfied with a large enough value given the noise level, and has been adapted from \cite{shi2022noise}.  Quasi-Newton Hessian approximations were preferred over the simpler choice of $H_k = \kappa_k I, \kappa_k > 0$ (which avoids the extra issues pointed above) because in practice this simpler alternative caused too slow convergence even without noise.
 
We test the average performance of the algorithm on a subset of small-scale CUTEr problems with general inequality constraints that contain less than 100 variables and less than 100 constraints, excluding the problems that contain equality or range constraints.  

We run the algorithm for a maximum of 1000 iterations by setting $\theta_1=0.1$, $\theta_2=0.01$.  The parameter $\Delta_k$ is set to a fixed value of $10^3$ for simplicity.  

We first run the algorithm without noise to get reference solution values by terminating it when $\|d_{qp}\|\leq 10^{-8}$.  Except the three problems for which the number of QN skips due to negative curvature was relatively large, the algorithm could provide a locally optimal solution in less than 50 iterations in the no-noise case.

\begin{figure}[h]
\scalebox{0.4}{\includegraphics{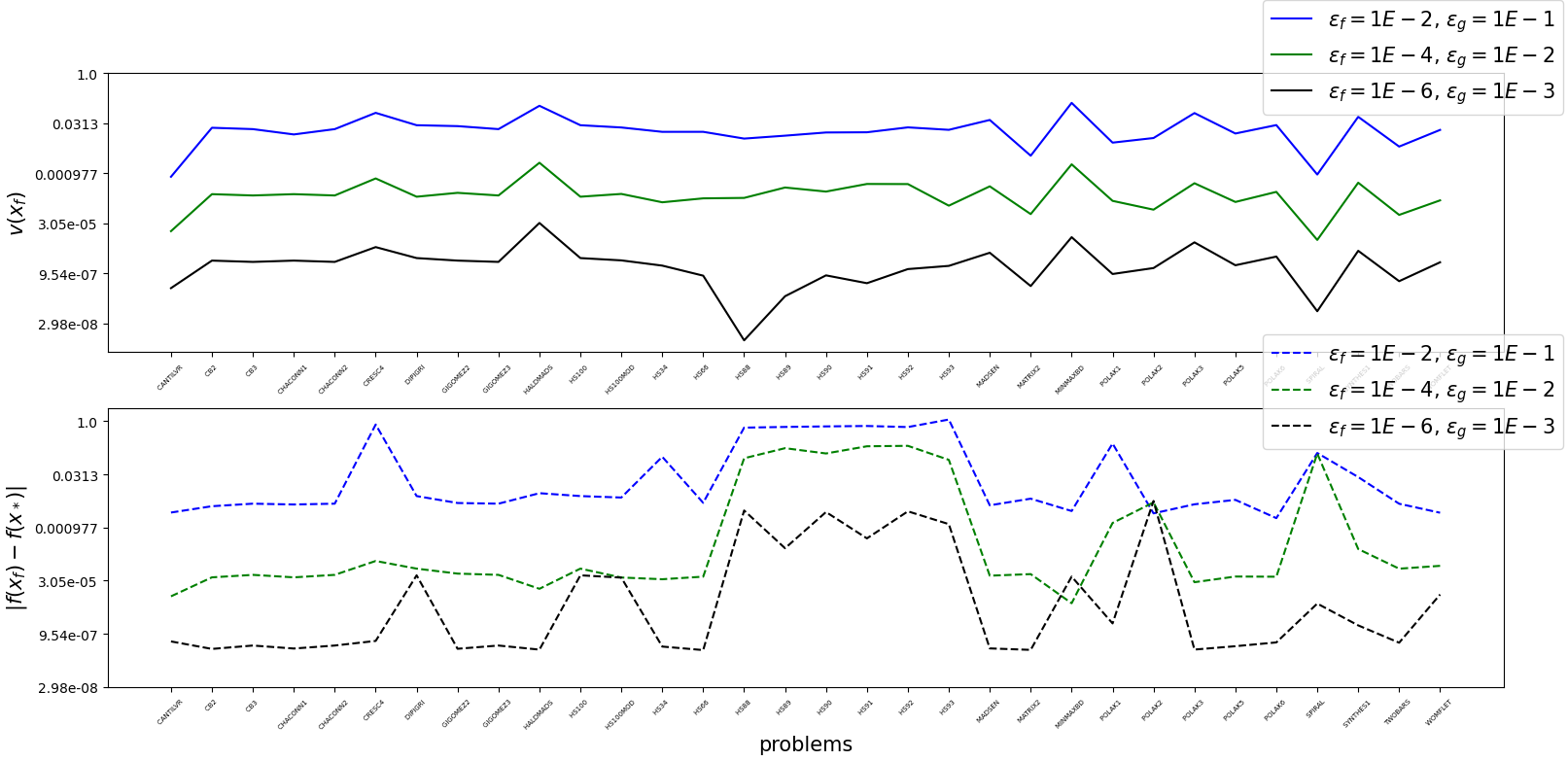}}
\caption{The dependence of $v(x_f)$ (upper subplot) and $|f(x_f)-f(x_\ast)|$ (lower subplot) to the noise level, where $x_f$ is the iterate returned by the algorithm.  The plots are generated at the logarithmic scale for clarity; the labels display the true (unscaled) values.}
\label{fig:averages}
\end{figure}

We next run the algorithm by injecting noise to function and gradient evaluations.  We added uniformly distributed random noise to each function evaluation (i.e. objective and constraint function evaluations), as well as each component of the objective gradient and constraint Jacobian.  That is, for function evaluations noise is a random variable with distribution $\mathcal{U}(-\epsilon_1, \epsilon_1)$, and for first order partial derivative evaluations noise is a random variable with distribution $\mathcal{U}(-\epsilon_2, \epsilon_2)$.  To get consistent results with the case of finite difference gradient approximations, we always set $\epsilon_2 = \sqrt{\epsilon_1}$ in our experimentation. 

In the presence of noise, it is not clear how to choose the final solution to be returned by the algorithm.  In our implementation, we set the point with the best (noisy) feasibility error in the last 100 iterations as the final solution.

We repeated our experiments three times with different realizations of the noise components.  Figure~\ref{fig:averages} summarizes the results.  From those plots, we can clearly observe the dependence of the solution quality to the noise level.  It is also interesting to observe the (generally) parallel alignment of the error lines, which to our understanding shows the effect of problem-specific parameters appearing in the theoretical bounds.   

\begin{figure}[h]
\scalebox{0.4}{\includegraphics{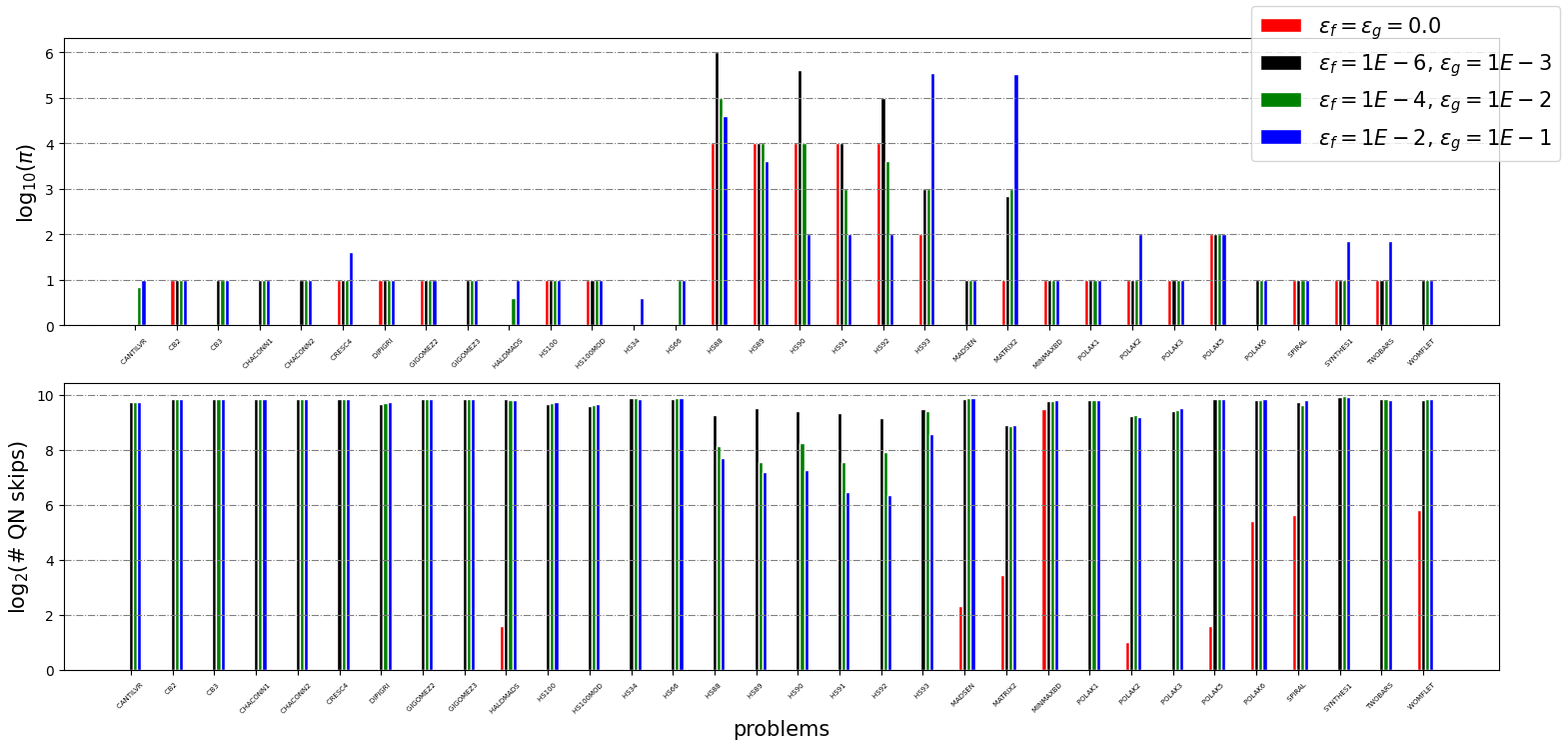}}
\caption{The final value of the penalty parameter $\pi$ (upper subplot) and the number of quasi-Newton update skips (lower subplot) for varying noise levels.}
\label{fig:more}
\end{figure}

We also investigated how noise effected the final value of the penalty parameter $\pi$.  We plot our findings in Figure~\ref{fig:more}.  Depending on the problem, the penalty parameter can take larger values for small ($\epsilon_f=10^{-6}$), or large ($\epsilon_f=10^{-2}$) noise levels.  Interestingly, for a few problems, when noise is present, the penalty parameter does not become as large as the no-noise case.  This might be related to the potential regularization effect of noise.  We present in the same figure the information on the number of times the quasi-Newton update has been skipped due to noise (or negative curvature).  We observe that the number of skips is close to the maximum number of iterations (i.e.1000) for most of the instances, although there were a few skips in the no-noise case.  This indicates, the quasi-Newton skips primarily occur in the critical region and are due to noise.  

\begin{figure}[h]
\scalebox{0.35}{\includegraphics{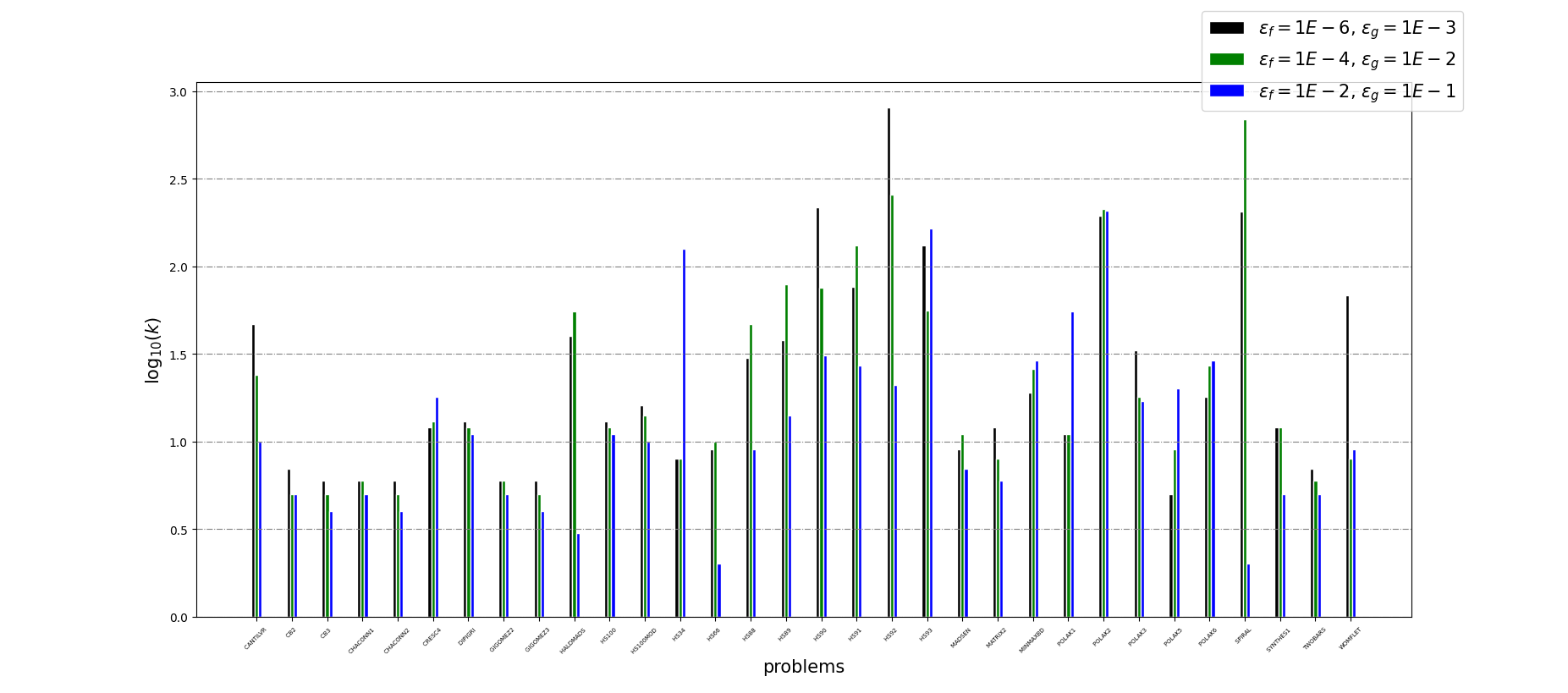}}
\caption{The iterate the algorithm finds for the first time a solution with an objective and constraint violation at most $2\epsilon_f$ away than $\tilde f(x_f)$ and $\tilde v(x_f)$.}
\label{fig:marked_iterates}
\end{figure}

To elaborate our observations on QN skips and to get an idea on how fast the algorithm arrives at a solution point comparable to the one reported after 1000 iterations, we check the iteration number at which the algorithm finds for the first time a solution with an objective value and feasibility error at most $2\epsilon_f$ away from the final values.  That is, we record
\[
k = \min \ k : \tilde v(x_k)\leq \tilde v(x_f) + 2\epsilon_f \quad \mbox{and} \quad \tilde f(x_k)\leq \tilde f(x_f) + 2\epsilon_f.
\]
Our observations are summarized in Figure~\ref{fig:marked_iterates}.  That iterate is rarely larger than 100, giving us an idea on how fast the algorithm arrives at the critical regions feasibility and optimality.  Although we run the algorithm for 1000 iterations, for most of the later iterations the change in the quality of the solution is in the order of the noise level.

\section{Conclusions}
In this paper, we studied a noise-aware SQP algorithm for optimization problems with noisy inequality constraints.  This work builds upon our earlier work on optimization with noisy equality constraints\cite{conDFO}.  The analysis present in this work is more extensive as it does not always assume satisfaction of regularity requirements.  The convergence of the algorithm has two components: a critical region for constraint violation and a critical region for optimality. To summarize in terms of $\epsilon=\max\{\epsilon_f, \epsilon_c, \epsilon_g, \epsilon_J\}$, we show that the iterates of the algorithm keep visiting the former region, where the stationarity measure for the minimization of constraint violation is $O(\sqrt{\epsilon})$.  Via similar lines of analysis, we show that if regularity holds this feasibility error bound improves to $O(\epsilon)$, and additionally the iterates infinitely often enter a critical region for optimality, where the optimality error is $O(\epsilon)$.  

On the practical side, our preliminary testing verifies the usefulness of noise-aware QN updates.  The only basic component missing for a practical implementation of the proposed algorithm is proper termination criteria.

\paragraph{Acknowledgments.} The authors are grateful to Jorge Nocedal for his helpful comments on an earlier version of this work.

\bibliographystyle{plain}

\bibliography{references}

\end{document}